\documentclass[twoside,12pt]{report}
\usepackage{amssymb,amstext,epsf}
\def\nc{\newcommand}

\nc\mpar[1]{}
\def\operatorname#1{\mbox{#1\,}}
\newtheorem{thm}{Theorem}[section] 
\newtheorem{hypothesis}[thm]{Hypothesis}
\newtheorem{lem}[thm]{Lemma}
\newtheorem{prop}[thm]{Proposition} 
\newtheorem{cor}[thm]{Corollary}
\newtheorem{fact}[thm]{Fact}

\newtheorem{defn}[thm]{Definition}
\newtheorem{remark}[thm]{Remark}
\newsavebox{\refs}
\def\extra{}
\def\newnum{}

\def\setproof#1#2{\def\extra{#1}\def\newnum{\sbox{\refs}{#2}}}
\newenvironment{proof}%
{\noindent {\proofstyle Proof\extra: } \sbox{\refs}{\record}%
 \global\def\extra{}\newnum\global\def\newnum{}%
 }%
{\hfill \fbox{\usebox{\refs}}}

\def\record{\arabic{chapter}.\arabic{section}.\arabic{thm}}

\def\proofstyle{\em \bf}

\def\bt{\begin{thm}} 
\def\et{\end{thm}}
\def\bl{\begin{lem}}
\def\el{\end{lem}} 
\def\bp{\begin{prop}} 
\def\ep{\end{prop}} 

\newcommand\fanc[1]{\mathbb{#1}}
 
\def\Z{\fanc{Z}} 
\def\reals{\fanc R} 
\def\cx{\fanc C}
\def\nats{\fanc N}
\def\integers{\fanc Z}
\def\zplus{\Z^+} 

\def\>{\rightarrow} 
\def\closure#1{\overline{#1}} 
\def\Int{\mathop{\text{Int}}}
\def\union{\cup} 
\def\nullset{\emptyset} 
\def\on#1{|_{#1}} 
\long\def\comment#1{}
\def\a{{\bf{A}}} 
\def\P{\{f^n(0)|n \in \zplus\}} 

\def\norm#1{\Vert #1 \Vert}
\long\def\ignore#1{}
\def\blackboard#1{{\bf{#1}}}
\def\term#1{\emph{#1}}
\def\crit#1{}
\def\newcrit#1{}

\def\picture #1 by #2 (#3){
  \vbox to #2{
    \hrule width #1 height 0pt depth 0pt
    \vfill
    \special{picture #3} 
    }
  }

\def\scaledpicture #1 by #2 (#3 scaled #4){{
  \dimen0=#1 \dimen1=#2
  \divide\dimen0 by 1000 \multiply\dimen0 by #4
  \divide\dimen1 by 1000 \multiply\dimen1 by #4
  \picture \dimen0 by \dimen1 (#3 scaled #4)}
  }

\author{Jeremy Kahn\thanks{Initial work supported by an DARPA NDSEG Fellowship}} 
\title{Holomorphic Removability of Julia Sets} 

\def\dist#1#2{{\cal QD}(#1,#2)}

\def\piece#1#2{P_{#2}(#1)}		
\def\annulus#1#2{A_{#2}(#1)}		
\def\taufunc#1#2{\tau_{#2}(#1)}  	
\def\critpiece#1{\piece{#1}0}
\def\critannulus#1{\annulus{#1}0}
\def\piecelevel{p}	

\def\statetiling{There exists an $L \in \zplus$ such that given any piece 
$P$ of
level greater than $L$, we can write 
$$
P = T \union R \union \bigcup (\closure{Q_i} \cap P),
$$ where $T$, $R$, and $\bigcup (\closure{Q_i} \cap P)$ are mutually
disjoint; $T$ is open, and $T \cap J = \nullset$; $R$ is compact and
holomorphically removable; and each of the $Q_i$ is a Yoccoz piece of
level $q_i>L$, the $Q_i$ are all mutually disjoint, and
$$f^{q_i-L}\on{Q_i}$$ is
univalent.  }

\def\lt{Lemma \protect\ref{tiling}}

\def\mt{{\cal C}}

\usepackage{fullpage}

\def\authorname{Jeremy Kahn}
\usepackage{fixtitle}
\date{\relax}

\begin{document}

\maketitle
\centerline{Abstract}
\bigskip

Let $f(z) = z^2 + c$ be a quadratic polynomial, with $ c$
in the Mandelbrot set $M$. 
Assume further that both
 fixed points of $f$ are repelling, and that $f$ is not
 renormalizable.  Then we prove that the Julia set $J_f$ of $f$ is
 holomorphically removable in the sense that every homeomorphism of
 the complex plane to itself that is conformal off of $J_f$ is in fact
 conformal on the entire complex plane.  As a corollary, we deduce
 that $M$ is locally connected at such $c$. 
\newcrit{Which was previously obtained by Yoccoz
\cite{Yoccoz:local:connectivity}}


\thispagestyle{empty}
\def\SBIMSMark#1#2#3{
 \font\SBF=cmss10 at 10 true pt
 \font\SBI=cmssi10 at 10 true pt
 \setbox0=\hbox{\SBF Stony Brook IMS Preprint \##1}
 \setbox2=\hbox to \wd0{\hfil \SBI #2}
 \setbox4=\hbox to \wd0{\hfil \SBI #3}
 \setbox6=\hbox to \wd0{\hss
             \vbox{\hsize=\wd0 \parskip=0pt \baselineskip=10 true pt
                   \copy0 \break%
                   \copy2 \break%
                   \copy4 \break}}
 \dimen0=\ht6   \advance\dimen0 by \vsize \advance\dimen0 by 8 true pt
                \advance\dimen0 by -\pagetotal
 \dimen2=\hsize \advance\dimen2 by .25 true in
%
%
  \openin2=publishd.tex
  \ifeof2\setbox0=\hbox to 0pt{}
  \else 
     \setbox0=\hbox to 3.1 true in{
                \vbox to \ht6{\hsize=3 true in \parskip=0pt  \noindent  
                {\SBI Published in modified form:}\hfil\break
                \input publishd.tex 
                \vfill}}
  \fi
  \closein2
  \ht0=0pt \dp0=0pt
 \ht6=0pt \dp6=0pt
 \setbox8=\vbox to \dimen0{\vfill \hbox to \dimen2{\copy0 \hss \copy6}}
 \ht8=0pt \dp8=0pt \wd8=0pt
 \copy8
 \message{*** Stony Brook IMS Preprint #1, #2 ***}
}

\SBIMSMark{1998/11}{December 1998}{}

\tableofcontents

\chapter{Statement of Main Theorem and Breakdown of Proof}
\label{overview}
\hyphenation{homeo-morph-ism}
\section{Introduction} \label{intro}
 Let $f(z)=z^2+c$ be a quadratic polynomial, with $c \in  M$
(where $M$ is the Mandelbrot set, defined in 
section \ref{Yoccoz Partition}). 
We consider two possible additional
hypotheses on $f$:

\begin{enumerate}
\item
Both of the fixed points of $f$ are repelling, and $f$ is not 
renormalizable;
\item
All of the periodic cycles of $f$ are repelling, and $f$ is not
infinitely renormalizable.
\end{enumerate}

\comment{
{\em The first hypothesis will be assumed throughout Chapters 1-3 \/}
 (excluding this section
(\ref{intro})); necessary modifications to the arguments for them to
 work with the second hypothesis will be described in Chapter
 \ref{further}.
}

Under either of the two above hypotheses, there are the following 
theorems:

\bt[Yoccoz]\label{jlc} $J_f$ is locally connected. \et

\bt[Yoccoz]\label{mlc} $M$ is locally connected at $c$. \et

\bt[Lyubich; Shishikura]\label{mzero} $J_f$ has measure 0. \et

See \cite{Yoccoz:local:connectivity}, \cite{Lyubich:area:zero}. See
also Milnor \cite{Milnor:local:connectivity} and Hubbard
\cite{Hubbard:local:connectivity} for expositions of Theorem \ref{jlc} 
(and also Theorem \ref{mlc} in the latter reference). 

\begin{defn}
We say that a compact subset $J$ of $\cx$ is {\em holomorphically
removable (HR)} in an open neighborhood $U$ of $J$ if, for every
topological embedding $h: U
\> \cx$, if $h\on{U-J}$ is conformal, then in fact $h\on
U$ is conformal.
\end{defn}
\begin{fact} \label{HR indep of K}
For each $K \ge 1$, $J \subset U$ is holomorphically removable if and 
only if
$J$ is removable for $K$-quasiconformal mappings, that is, 
 for every
topological embedding $h: U
\> \cx$, if $h\on{U-J}$ is $K$-quasiconformal, then in fact $h\on
U$ is $K$-quasiconformal.
\end{fact}
For a proof, see section V.3 of  \cite{Lehto:Virtanen}, 
where the conditions given here on $h$  are just those to put it in 
Lehto and Virtanen's 
class ${\cal W}_2$ of functions.

Clearly, if $J \subset U \subset V$, and $J$ is holomorphically
removable in $U$, then it is holomorphically removable in $V$. Using
Fact \ref{HR indep of K} above, it is easy to show that the converse
is true, that $J$ is HR in $U$ if it is HR in $V$ (assuming of course
that $J$ is compact).  Thus we can suppress mention of the
neighborhood and just assume $U=\cx$. 

The simplest example of a
holomorphically removable set is a point. The next simplest is a
piecewise smooth curve.

The purpose of this work is to prove the following theorem (with the
same hypotheses): 
\bt[Main Theorem] \label{HR} 
$J_f$ is holomorphically
removable. 
\et
In section \ref{MLC} we give use Theorem \ref{HR} to give a quick proof of
Theorem \ref{mlc}. {\bf Throughout the first three chapters will we always
assume the first hypothesis on $f$.} The proof with the weaker second hypothesis
will be discussed in section \ref{fr}. We mention that the critically
non-recurrent cases (see chapter
\ref{decomp chapter} for a definition) also follow from the work of Jones
\cite{Jones:hr,Jones:Carleson:Yoccoz}.
Speculations on further holomorphic removability results are discussed in 
section \ref{conj}.

\subsection{Acknowledgements}
I would like to thank Curtis McMullen, Mikhail Lyubich, and Mitsuhira Shishikura
for helpful conversations.
\newcrit{Curt showed me Yoccoz's trick (or at least hinted at it), Lyubich
provided encouragement early on, and Shishikura explained the connection to his
argument, thus providing a new and different argument for the tiling lemma.}

\ignore{
\section{Overview of proof}

We now offer for the reader an overview of the entire proof
the main Theorem, \ref{HR}, assuming hypothesis 1.

Let $J=J_f$ be the Julia set of the given quadratic poymonial $f$. 
Then given a homeomorphism $h:\cx \> \cx$, with $h \on {\cx-J}$ 
conformal,
we need to show that $h$ is conformal on all of $\cx$.
We prove this  by uniformly approximating such a homeomorphism $h$ by
uniformly quasiconformal mappings.  It then follows that the
homeomorphism is quasiconformal; a more careful examination shows in
fact that it must be conformal.

These uniformly quasiconformal approximations are obtained via
Yoccoz's Markov-like partition of the dynamical plane. We define
conformally invariant ``quasiconformal distortion bounds'' for a pair
consisiting of a Jordan domain and a closed subset of it: these
distortion bounds control the dilatation of a quasiconformal
mapping with the same boundary values as an arbitrary re-embedding of
Clthe Jordan domain that is conformal on the complement of the closed
subset.

Our main proposition states that the quasiconformal distortion bounds
for the pair consisting of a Yoccoz partition piece and the
intersection of the piece with the Julia set is uniformly bounded above,
independent of the piece and depending only on $f$. This proposition
provides the required $K$-quasiconformal approximations to the given
homeomorphism, with $K=K(f)$ depending only on $f$, and not on
the homeomorphism. These approximations show that the homeomorphism is always
$K$-quasiconformal. This in turn implies $J_f$ has zero 2-dimensional
Lebesgue measure, and therefore the homeomorphism is conformal.

Two lemmas are used in the proof of this main proposition of uniform
distortion bounds.  One is a piece-dependent version, where the
distortion bounds are allowed to depend on a given piece.  The other
is a tiling lemma, which describes how every piece can be ``tiled'' by
univalent dynamical copies of elements of a fixed finite set of
pieces.  These copies do not necessarily tile (cover) the whole
piece, or even cover the whole intersection of the Julia set with the
given piece, but the leftover portion of the Julia set is itself shown
to be holomorphically removable. Then, every piece is shown to have qc
distortion bounds no greater than the supremum of the bounds of a
fixed finite set of pieces, as follows: given a re-embedding of the
given piece that is conformal off of the piece's intersection with the
Julia set, we can replace the re-embedding on each of the tiling pieces 
with a
$K$-quasiconformal map with the same boundary values on that ``tile'',
where $K$ is the supremum of the distortion bounds of the fixed finite
set of pieces. The resulting map is then $K$-quasiconformal on the
complement of the leftover untiled portion of the Julia set. But this
set is holomorphically removable, and so the resulting map is
$K$-quasiconformal on the whole piece, as desired.

The proof of the main theorem is thus reduced to two lemmas: the
non-uniform piece-dependent distortion bounds, and the tiling lemma.

To obtain piece-dependent distortion bounds, we quasiconformally embed
copies of a canonical Jordan domain and closed-subset pair into the
given piece, so that the images of the closed subsets cover all of
the Julia set in that piece. We
prove qc distortion bounds for this canonical model.  The distortion
bounds for the given piece then follow.

To prove the tiling lemma, we use the Yoccoz puzzle theory. We examine
annuli that are differences between two nested puzzle pieces of
consecutive levels. The Yoccoz theory then tells us that the sum of
the moduli of such annuli that surround the critical point diverges.
Annuli are pulled back to form a disjoint collection of annuli, with
divergent total modulus around every point of the leftover untiled
portion of the Julia set. (The Yoccoz theory, as refined by
Shishikura (and also by Lyubich), also describes how to pull back these
annuli). This implies that this leftover portion is holomorphically
removable.
}

\section{Yoccoz Partition} 
\label{Yoccoz Partition}
The main tool for proving all the above theorems is the
Yoccoz partition, which we now describe
 \cite{Milnor:local:connectivity, Hubbard:local:connectivity}.
  Let us
 first recall some basic theory and terminology for the dynamics of
 quadratic polynomials. 
Given a quadratic polynomial $f(z)= f_c(z)= z^2 +c$, let 
$K(f)= \{z \mid f^n(z) \not\> \infty\}$. 
Then $J_c=J(f)= \partial K(f)$, 
and $K(f)$ is connected if and only if $0 \in K(f)$. 
Under hypothesis 1 (or 2) on $f$, $K(f)=J(f)$.
Then the Mandelbrot set $M$ is defined by 
$$
M=\{c \mid 0 \in K(f_c) \}. 
$$  
If $K(f)$ is connected, then there exists a unique conformal isomorphism
$\phi : \cx - \closure {\Delta} \> \cx - K(f)$ for which 
$\phi (z^2) = (\phi(z))^2 +c$.  
An \term{external ray} $R(\theta)$ is then defined by 
$$
 R(\theta) = \phi(\{r e^{2 \pi i \theta} \mid 1 < r < \infty\}).
$$  
The map $f$ acts as angle doubling modulo 1 on the external rays 
$R(\theta)$:
$f(R(\theta))=R({2\theta})$. 
We say that $R(\theta$) \term{lands} at $z \in J_c$ if 
$\lim_{r \> 1} \phi(r e^{2 \pi i \theta}) = z$. 
We first recall \cite{Milnor:dynamics:lectures,
 Hubbard:local:connectivity} two basic results about the landing of
 external rays:
\bp 
If $\theta$ is periodic under doubling modulo 1, then $R(\theta)$ lands 
at a
parabolic or repelling periodic cycle. Conversely, if $z \in J$ is a
repelling (or parabolic) periodic point, then at least one periodic
ray lands at it. In the case where $z$ is a fixed point, then set of rays landing
at $z$ are cyclically  permuted by the $f$.
\ep
  We are assuming both fixed points of $f$ are repelling. The
zero external ray lands at one of them, called the $\beta$ fixed
point, or just $\beta$. The other fixed point is called  $\alpha$. At
least two rays land at $\alpha$ (because the only cycle of length 1
that is periodic under doubling is $\{0\}$, which lands at
$\beta$). Form the connected 1-complex $\Gamma_0$ consisting of
$\alpha$, the portion of the rays landing at $\alpha$ with potential
less than 1 (the \term{potential} of a
point $z \in \cx -K_f$ is defined as $\log |\phi^{-1}(z)|$), and the
equipotential curve of potential 1. For $n \in
\zplus$, let $\Gamma_n=f^{-n}(\Gamma_0)$. A {\em piece} of level $n$
is a bounded component of $\cx-\Gamma_n$. Each piece is a Jordan
domain.  If $n<m$, then $f^{m-n}$ maps every piece of level $m$ to a
piece of level $n$.

In Yoccoz's work \cite{Milnor:local:connectivity,
Hubbard:local:connectivity}, Theorem \ref{jlc} is proven by showing
the following, which will be used in Chapter \ref{decomp chapter}:
\bt \label{dzero} 
The diameter of all pieces of level
$n$ goes uniformly to zero as $n \> \infty$. 
\et 

\setproof{ of \ref{jlc}}{\ref{jlc}}
\begin{proof}
To show that $J$ is locally connected at a given point $z \in J$ (with
$f^n(z) \neq \alpha$ for all $n$), consider the pieces of all levels
that contain that point.  They are connected and open, and
\ref{dzero} above tells us that that they form a neighborhood base for
$z$.  (In the case where $f^n(z)=0$ (so $z \in \Gamma_n$) for some
$n$, the interior of the union of closures of pieces of level $n$ that
border on $z$ is a connected open neighborhood of $z$, with diameter
going to zero as $n \> \infty$).
\end{proof}

The theory used to prove \ref{dzero} will be discussed in 
Chapter \ref{decomp chapter}, where it will be used to show further 
results.

\section{Quasiconformal Distortion Bounds} 
\label{distortion bounds}

Let us now introduce the general concept of quasiconformal distortion
bounds. Let $U \subset \cx$ be a Jordan domain, and $A$
a closed subset of $U$. Suppose there exists $K$ such that for all
embeddings $h: \closure{U} \> \cx$ with $h\on{U-A}$ conformal, there
exists an embedding $\tilde{h}:\closure{U} \> \cx$ such that
$\tilde{h}\on U$ is $K$-qc, and $h\on{\partial
U}=\tilde{h}\on{\partial U}$. Then we let $\dist A U$ be the
least such $K$ (and set $\dist A U = \infty$ if there is no such $K$).
In practice we will just be interested in establishing upper bounds
for $\dist A U$, or just showing that it is finite. We call such
bounds {\em qc distortion bounds}.

For future reference, we include some basic facts about these distortion 
bounds:

\begin{fact} \label{conf inv}
$\dist A U$ is a conformal invariant: if $g: U \> V$ is a conformal
isomorphism with $g(A)=B$, then $\dist A U = \dist B V$
\end{fact}
\begin{proof}
By Caratheodory's theorem\cite{Milnor:dynamics:lectures}, $g$ extends to
homeomophism between
$\closure U$ and
$\closure V$. The result then follows immediately. 
\end{proof}
\begin{fact} \label{monotonicity}
 If $A \subset B \subset U$, then $\dist A U \le \dist B U$.
\end{fact}
This is immediate.

The following fact shows that it can  be sufficient to assume that $h$
is only quasiconformal:
\begin{fact} \label{replqc}
Suppose $\dist A U \le K$, and $h: \closure{U} \> \cx$ is an embedding
with $h\on{U-A}$ $L$-qc. Then there exists an embedding
$\tilde{h}:\closure{U} \> \cx$ such that $\tilde{h}\on U$ is $KL$-qc,
and $h\on{\partial P}=\tilde{h}\on{\partial P}$.
\end{fact}


\begin{proof}
 Let the Beltrami coefficient $\mu$ on $h(U)$ be equal to the complex
dilatation of $h^{-1}$ on $h(U-A)$, and zero on $A$.  Let $g: h(U)
\> V$ be a quasiconformal map with complex dilatation $\mu$.  (The
existence of $g$ is guaranteed by the Measurable Riemann
Mapping Theorem
\cite{Ahlfors:Bers}. We can assume that $V$ is a Jordan domain, and that $g$
extends to a homeomorphism $g: h(\closure U) \> \closure V$).
  Then $g$ is $L$-quasiconformal, and $g \circ h: \closure U \>\closure V$ is a
homeomorphism that is conformal on $U-A$.  Therefore there exists
$\widetilde{(g \circ h)}: U \> V$ that is $K$-quasiconformal, and
agrees with ${g \circ h}$ on $\partial U$.  So let $\tilde{h} = g^{-1}
\circ \widetilde{(g \circ h)}$: it has the required properties.
\end{proof}

\begin{fact} \label{transfer}
If there is a homeomorphism $g: \closure U \> \closure V$ such that $g \on U$
is $L$-qc, and $g(A)=B$, then $\dist B V \le L^2 \dist A U$.
\end{fact}
\begin{proof}
Given $h: \closure V \> \cx$ with $h \on {V-B}$ conformal, let
$\tilde h= \widetilde{(h \circ g)} \circ g^{-1}$. Here $\widetilde{(h 
\circ g)}$
is as given from $h \circ g$ by Fact \ref{replqc}.
\end{proof}

We can also state and prove a more general fact:
\begin{fact} \label{clever transfer}
If there is a homeomorphism $g: \closure U \> \closure V$ with $g(A)=B$ and
$g \on {U-A}$ $L$-quasiconformal,
 then $\dist B V \le L^2 (\dist A U)^2$.
\end{fact}
\begin{proof}
Given an embedding $h: V \> \cx$ with $h \on {V-B}$ conformal, we must
find a $L^2 (\dist A U)^2$-qc map $\tilde h: V \> \cx$ with $h \on
{\partial V}= \tilde h \on {\partial V}$.  Now $h \circ g: U \> \cx$
is an embedding that is $L$-qc on $U-A$, so by Fact \ref{replqc} there
exists a $L \cdot \dist AU$-qc map $\widetilde{(h \circ g)}: U \> \cx$
that agrees with $h \circ g$ on $\partial U$.

Now note that, by Fact \ref{replqc}, there exists a $L\cdot\dist A U$-qc 
map
$\tilde g: \closure U \> \closure V$ with $\tilde g \on {\partial U} = g \on {\partial U}$.

Then $\widetilde{(h \circ g)} \circ \tilde g^{-1}: \closure U\> \cx$ is $L^2
(\dist A U)^2$-qc (on $U$), and agrees with $h$ on $\partial U$. It is the 
required map
$\tilde h$.
\end{proof}

\begin{fact} \label{compactly}
If $A \subset U$ is compact, then $\dist A U \le \infty$.
\end{fact}
\begin{proof}
By the Riemann mapping theorem (and Caratheodory's theorem), we can assume $U$
and $h(U)$ are both the unit disk, and $0\in h(A)$. Then, using Schwartz
reflection, we find that $h\on{S^1}$ is real-analytic, and, using Montel's
theorem (or the Koebe distortion theorem), that $h' \on {S^1}$ is
bounded.  Likewise for $h^{-1}\on {S^1}$ (one checks that $h(A)$ always 
lies
within some definite subdisk (depending only on $A$), because $h(U-A)$ 
has some
fixed modulus). Therefore the map on the boundary is uniformly 
bi-Lipschitz,
which is certainly enough to insure a uniformly quasiconformal extension
(e.g. just cone it off, mapping $(r, \theta)$ to $(r, h(\theta))$ ).
\end{proof}

Fact \ref{conf inv} will be used in section
\ref{two together}; the others  will be used in Chapter
\ref{non-uniform chapter}.

\section{Uniform Distortion Bounds}\label{udb}
The proof of the Main Theorem, \ref{HR}, can be reduced to the
following lemma:
\bl[Uniform Qc Bounds] \label{uniform} 
There exists a $K$, depending only on $f$, such that
for all pieces $P$, $\dist {J \cap P} P \le K$.  \el

Assuming this Lemma, we can complete the proof of Theorem \ref{HR}:

\setproof{}{\ref{HR}}
\begin{proof}
We will first show that there exists a $K$  such that if $h:\cx \> \cx$ is a
homeomorphism, with $h\on{\cx-J}$ conformal, then $h$ is $K$-quasiconformal.
This $K$ will be independent of $h$. 

Let $K$ be as given in Lemma  \ref{uniform}. We will show that $h$ is
$K$-quasiconformal by approximating it  uniformly with $K$-qc maps.
For each $n \in \zplus$, we define $h_n: \cx\>\cx$ as follows:
let $h_n = h$ on the unbounded component
of $\cx -
\Gamma_n$, and for each piece $P_i$ of level $n$, let $h_n = 
\tilde{h}_i$ on
$\closure{P_i}$, where $\tilde{h}_i\on{\partial P_i}=h\on{\partial
P_i}$, and $\tilde{h}_i\on{P_i}$ is $K$-qc. (The existence of the
$\tilde{h}_i$'s are guaranteed by Lemma
\ref{uniform}). Then $h_n$ is $K$-qc. (Here we use the fact that 
$\Gamma_n$, a
piecewise smooth 1-complex, is holomorphically removable). Now,
because the diameters of the pieces goes to zero as $n \> \infty$, so
do their images by $h$. But $h_n$ maps every piece of level $n$ to its
image under $h$. Therefore $\Vert h-h_n
\Vert_\infty \> 0 $ as $n \> \infty$. So $h$ is $K$-qc.   

Now, the above fact (there exists $K$ such that given $h:\cx \> \cx$ a
homeomorphism, $h\on{\cx-J}$ conformal, then $h$ is
$K$-quasiconformal) implies that $J$ has zero area (thus we have also
proven Theorem \ref{mzero}). For if not, one can take any Beltrami
coefficient supported on J with dilatation (that is, essential
supremum of pointwise dilatation) greater than $K$, and using the
Measurable Riemann Mapping Theorem \cite{Ahlfors:Bers}, integrate it
to obtain a quasiconformal homeomorphism of $\cx$ that is conformal
off of $J$ but has dilatation greater than $K$, a contradiction. 
We have thus shown so far that any homeomorphism $h:\cx \> \cx$ with
$h\on{\cx-J}$ conformal is ($K$-)quasiconformal, and conformal off of
a set of measure 0. But then, we can conclude, as wanted, that it is
conformal, by the following\cite{Ahlfors:book:qc, Lehto:Virtanen}:
\bt A quasiconformal mapping that is conformal off of a set of measure 
zero is
conformal. \et
\end{proof}

\section{Proof of the Uniform Distortion Bounds}
\label{two together}

 There are two lemmas that form the basis of the proof of Lemma
$\ref{uniform}$. (One will be proven in each of the following two
chapters). The first is a non-uniform version, where we allow $K$ to
depend on the piece:
\bl[Piece-dependent Qc Bounds] \label{nonuniform}
 For all pieces $P$, there exists a $K(P)$ such that $\dist {J \cap P}
P \le K(P)$.
\el
The second lemma breaks down each piece into copies of pieces at a 
fixed
level:

\bl[Tiling Lemma] \label{tiling}
\statetiling
\el
\begin{remark}
We allow either a finite or countable set of $Q_i$'s, typically the 
latter.
\end{remark}
\begin{remark}
Thus each $Q_i$ is a univalent copy of a piece at level $L$, by a map
(namely, an iterate of $f$) that maps Julia set to Julia set.
\end{remark}

\setproof{ of Lemma \ref{uniform}, given Lemmas \ref{nonuniform} and
\ref{tiling}}{\ref{uniform}}
\begin{proof}
Let $L$ be as given by \lt. Then let $K$ in Lemma \ref{uniform} to
be the maximum of the $K(P)$'s of the (finitely many) pieces of level
at most $L$ (as given by Lemma
\ref{nonuniform}). We will show that Lemma
\ref{uniform} holds for this choice of $K$.

This $K$ works tautologically for all pieces of level at most $L$. Now
let $P$ be a piece of level greater than $L$. By \lt, we may write
$$
P = T \union R \union \bigcup (\closure{Q_i} \cap P).
$$  
Then, given $h: \closure{P} \> \cx$, we define $\tilde{h}$ as
follows:

For each $Q_i$,
\begin{eqnarray*}
\dist{J \cap Q_i}{ Q_i} &=& \dist{f^{q_i-L}(J \cap Q_i)}{ 
f^{q_i-L}(Q_i)}\\
 \lefteqn{\text{(because $f^{q_i-L} \on {Q_i}$ is univalent)}}\\ 
 &=& \dist{J \cap f^{q_i-L}(Q_i)}{
f^{q_i-L}(Q_i)} \\
 &\le& K(f^{q_i-L}(Q_i)) \le K.
\end{eqnarray*}
So we can replace $h \on{ \closure{Q_i}}$ by $h_i$,
 with $h_i \on{ \partial Q_i} = h \on{ \partial Q_i}$, and $h_i \on{ 
Q_i}$
$K$-quasiconformal.

Define $\tilde{h}$ by $\tilde{h} \on{ \closure{Q_i}} = h_i \on{ 
\closure{Q_i} }$,
	 and $\tilde{h} = h$ off of $\bigcup \closure{Q_i}$. 
Then $\tilde h$ is well-defined and continuous on $\bigcup \closure{Q_i}$ 
(because $\tilde h = h$ on $\bigcup \partial Q_i$), and
$\tilde h$ is continuous on $\closure{\bigcup Q_i}$, since the
diameters of the $Q_i$ (and their images under $h$, $\tilde h$) goes
to zero as $i \rightarrow \infty$.  So $\tilde h$ is continuous on
$\closure P$. It is also injective, and hence is an embedding

We now just need to verify that $\tilde h$ is $K$-qc on $P$.  First
note that it is $K$-qc on $T \cup \bigcup Q_i$.  Therefore it is
$K$-qc on the open set $P-R=T \cup \bigcup (\closure{Q_i} \cap P)$,
because $\bigcup (\partial Q_i \cap P)$ is a piecewise smooth locally
finite 1-complex.  Therefore it is $K$-qc on $P$, because the
remaining set, $R$, is holomorphically removable.
\end{proof}

Chapters \ref{non-uniform chapter} and \ref{decomp chapter} give the
	proofs of Lemmas \ref{nonuniform} and \ref{tiling} repectively, thus
	completing the proof of Theorem \ref{HR}.

\chapter{The Piece-dependent Bounds}
\label{non-uniform chapter}

In this chapter we prove the piece-dependent distortion bounds.  We
introduce a canonical model, and prove
quasiconformal distortion bounds for it. Then, given an arbitrary
Yoccoz puzzle piece $P$, we embed this canonical model into $P$ in such a way as
to imply qc distortion bounds for $P$.

In section \ref{nonuniform bounds breakdown} we define this canonical model
and describe its role in the proof of piece-dependent distortion bounds. In
section \ref{bounds for rns} we prove qc distortion bounds for it. In
section \ref{embedding rns} we describe how it is embedded into a given
piece $P$. How all of this fits together to prove piece-dependent
distortion bounds for $P$ is also described in section \ref{nonuniform
bounds breakdown}.

\section{The Role of the ``Recursively Notched Square''}
\label{nonuniform bounds breakdown}

\begin{figure}
\centerline{\epsfbox{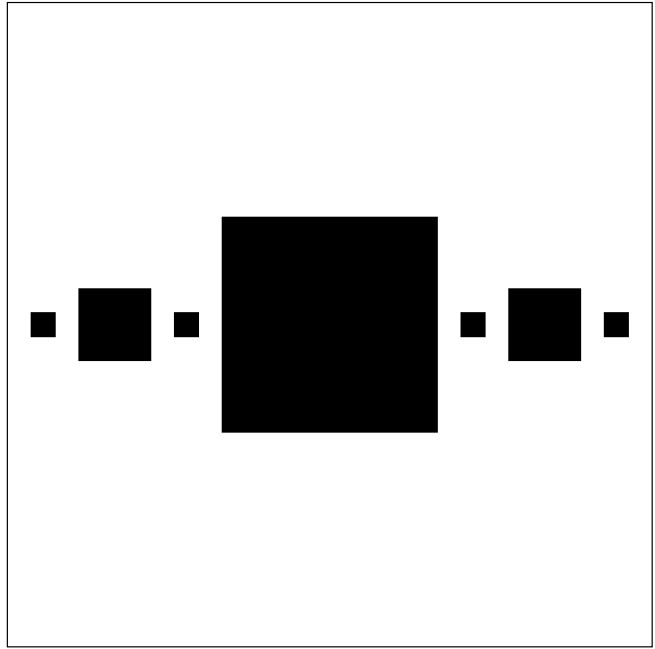}}
\caption{The recursively notched square}
\label{rns}
\end{figure}
\mpar{Make RNS into RNR ($[0,1] \times [-1,1]$)?}
We first define the ``recursively notched square'' as the pair
$(S,N)$, which are defined as follows. 
\crit{Curt will want the English changed here}
Take the open square $S=(0,1)\times(-1/2,1/2)$, and divide it into nine
equal-sized smaller squares in the obvious way. There are unique
homotheties (i.e.  direction-preserving similarities) from the large
square to each of the smaller squares. Define $N$ to be the smallest
subset of the $S$ such that $N$ contains the central small square, and
$N$ contains its own image under the homotheties $h_l$ to the
middle-left square and $h_r$ the middle-right one. 
(We have $h_l(z)=z/3$ and $h_r(z)=(z-1)/3 + 1$.) 
 So if we let
$(a_i)_{i=1}^n$ denote a sequence of $l$'s and $r$'s, then
$$
N = \bigcup 
h_{a_1}\circ h_{a_2} \circ \dots \circ h_{a_n} (\closure{S})
$$
where the union ranges over all sequences of length $n \ge 0$. 
See figure \ref{rns}. 
Note that $(\closure S - \Int N) \cap \reals = \mt$, 
\mpar{Need better notation for $M$}%
where $\mt \subset [0,1]$ denotes  the
middle-thirds Cantor set
$\{\sum_{i=1}^\infty a_i 3^{-i} \mid a_i \in \{0,2\} \}$.

The key step toward showing Lemma
\ref{nonuniform} is the following:

\bl \label{Sbound} 
The recursively notched square has quasiconformal distortion bounds:
$$
\dist {\closure N \newcrit{\closure N \cap S?}} S = D_0 < \infty.
$$ 
\el

Now\mpar{Already know from 1.3, right?}
if $P$ is any level $n$ Yoccoz piece,
$\partial P \cap J$ is a finite set because it is a subset of $f^{-n}(\alpha)$.  We
will show that we can cover a neighborhood of each point in this set by a copy
of $(S,\closure N)$\newcrit{$(\closure S,\closure N)$?}. More precisely,

\bl \label{getS}
 Given $P$, we can find $K, h_1, \dots, h_m$ and open subsets $R_1, \dots, R_m$ of
$P$  such that $h_i:\closure S \>\closure{ R_i}$ is a homeomorphism with $h_i \on
{S-\closure N}$ $K$-qc, $h_i(\closure N) \supset R_i \cap J$, the $\closure{R_i}$ are
disjoint (and have holomorphically removable boundary), and $(P
\setminus (\bigcup R_i)) \cap J$ is compactly contained in $P$. \el

Given Lemmas \ref{Sbound} and \ref{getS} we can now prove Lemma
\ref{nonuniform}, using the basic facts from section \ref{distortion bounds}:

\setproof{ of Lemma \ref{nonuniform}}{\ref{nonuniform}}
\begin{proof}
We have that $\dist{J \cap R_i} {R_i} < K^2 D_0^2$ by Facts \ref{clever transfer} and 
\ref{monotonicity}. Therefore,
given an embedding $h: \closure P \> \cx$ with 
$h \on {P-J}$ conformal, we can 
replace $h$ on each $R_i$ with a 
$K^2 D_0^2$-quasiconformal map with the same boundary values, 
and thereby obtain an embedding that agrees  with $h$ on $\partial P$, 
and which is $K^2 D_0^2$ quasiconformal off of a compact subset of $P$, namely 
$E= (P- \cup R_i) \cap J$. 
By Fact \ref{compactly}, $\dist E P$ is finite, say $K_1$. Then by
 Fact \ref{replqc}, there exists a $K_1K^2D_0^2$-quasiconformal map
 $\tilde h: \closure P \> \cx$ with $\tilde h \on {\partial P} = h \on
 {\partial P}$. So
$$
\dist {J \cap P} P \le K_1K^2D_0^2 < \infty.
$$
\end{proof}

\section{Qc Distortion Bounds for the RNS}
\label{bounds for rns}


\def\f{\xi} 
To prove quasiconformal distortion bounds for the recursively notched square,
we will
in fact show a stronger property, 
 from which quasiconformal distortion
bounds can be deduced.
\def\Sobolev{W^{1,2}}
\def\normS#1{\norm{#1}_{1,2}} Let $U \subset \cx$ be open. 
We denote by $\Sobolev(U)$
the Sobolev space of functions $\f: U \> \reals$ (modulo
constants) with one distributional derivative in $L^2$, with norm
$$
\norm{\f}^2_U = \norm{\f}_{U,1,2}^2= \int\!\int \left(\frac{\partial
\f}{\partial x}\right)^2 +
\left(\frac{\partial \f}{\partial y}\right)^2 \,dx\,dy =
2i\int\!\int \left(\frac{\partial \f}{\partial z}\right)
\left(\frac{\partial \f}{\partial \overline{z}}\right)\,dz\,d\overline{z}.
$$ 
\begin{remark}
The usual norm for $\Sobolev(U)$ also includes the usual $L^2$ norm, obviating the
need to mod out by constants. But what we need is the above.
\end{remark}
\begin{remark}
\mpar{reference? (ask Tom)}
Functions in this space are not necessarily continuous.\newcrit{``We will state the
necessary continuity properties later''?}
\end{remark}
The latter formula shows the norm is conformally invariant, in
the sense that if
$h: V \> U$ is conformal, then $\norm{\f\circ h}_V=\norm{\f}_U$. When there is
no danger of confusion, we will omit the domain in the norm notation. 
Furthermore, the norm is quasiconformally quasi-invariant:
\begin{fact} \label{q q}
If $h: V \> U$ is $K$-quasiconformal, 
then $\norm{\f\circ h}_V \le K \norm{\f}_U$.
\end{fact}
\begin{proof}
Suppose $\f$ is $C^1$ with compact support; then $\f \circ h \in \Sobolev$ 
because $h$ is in $\Sobolev$ and is absolutely continuous for 2-dimensional Lebesgue
measure (and thus the change-of-variable formula applies). An easy calculation
then verifies the inequality in this case (see \cite[Ch.~1,
Sec.~F]{Ahlfors:book:qc}). But such $\f$ are dense in $\Sobolev$, so the result
follows.
\end{proof}

\ignore{
This is shown in the case where $h$ is differentiable in the last section
of the first chapter of \cite{Ahlfors:book:qc}. It is used without concern
for differentiability issues in \cite{Nag:Sullivan:quantum}\newcrit{check this?}, in
fact in the proof there of Theorem \ref{D-bounded implies qc}. It is not hard to show
it in the general case, using the facts that $h$ itself has
one distributional derivative in $L^2$, and that $h$ is absolutely continuous 
for 2-dimensional Lebesgue measure. \newcrit{I should think about this, maybe write
it down}
}
\crit{Is it actually shown later in Ahlfors? I doubt it}

\def\distS#1#2{{\cal SD}(#1,#2)} Now, suppose again we are given  $A  \subset
U$ closed. We define $\distS A U$ as the least $K$ such that, for all
$\f: \closure{U} \> \cx$ continuous, with $\normS{\f \on {U-A}} \le 1$,
there exists $\tilde{\f}:\closure{U} \> \cx$ continuous such that 
$\tilde{\f} \on {\partial U}=\f \on {\partial U}$, and
 $\norm{\tilde{\f}
\on U} \le K$.

\begin{prop}
 For all $K$ there exists $K'$ such that for all $A \subset U
\subset \cx$, if $\distS A U \le K$,  then $\dist A U \le K'$.
\end{prop}

\begin{proof} We will use a result of  Nag and Sullivan\cite{Nag:Sullivan:quantum},
which states:
\begin{thm}\label{D-bounded implies qc}
Suppose that $X$ and $Y$ are Jordan domains in $\cx$, and $h:\partial X \>
\partial Y$ is an orientation-preserving homeomorphism.
Suppose there exists a $C$ such that for all $f$ continuous on
$\closure Y$ with $\norm{f \on Y}_{1,2} \le 1$, there exists
a continuous extension $g$ (to $\closure X$) of $\,f\on
{\partial Y} \circ h\,$ with $\norm{g \on X} \le C$.  Then
$h$ has an extension $\tilde h: \closure X \> \closure Y$ such that 
$\tilde h \on X$ is 
$C'$-quasiconformal, with $C'$
depending only on $C$.
\end{thm}

Now, given an embedding $h$ of $\closure U$ that is conformal on $U-A$, let
$V=h(U)$.  For all continuous functions $f$ on $\closure V$ with
$\norm{f}_{1,2} \le 1$, we find that $\norm{f \circ h \on {U-A}} \le
1$, and therefore we can find $g$ continuous on $\closure U$ with $g
\on {\partial U} = f\circ h \on {\partial U}$, and with $\norm{g \on U} \le
K$. Using the theorem above, we conclude that $\partial h$ has a
$K'$-quasiconformal extension, with $K'$ depending only on $K$.
\end{proof}

To prove quasiconformal distortion bounds for the recursively notched
square, we will show that $\distS {\closure N} S \le \infty$.

\def\strip{{\cal F}}
Now let $\strip=\{ z |0 < \Im z < \pi\}$ be an infinite strip. Let us
define $\hat{\strip}$ as $\strip\cup \partial \strip$, where $\partial
\strip$ denotes the ideal boundary of $S$. Then $\hat{\strip}$ may be
identified with the closure of $\strip$ in $\cx$, plus two points,
positive and negative (real) infinity, with the obvious neighborhood
bases.

\def\slits{V}
We will show:
\begin{lem}[Mapping lemma] \label{mapping}
There is a quasiconformal homeomorphism $h: S-N \> \strip -
\closure{\slits}$, where $\slits$ is a union of  countably many vertical slits in
$\strip$ such that
\begin{itemize}
\item the imaginary part of each of the slits is bounded between
$\pi/5$ and
$4\pi/5$, and
\item $\closure{\slits} \subset (\slits \cup M)$ (where $M=\{z |
\Im z = \pi/2\}$ is the midline of $\strip$). 
\end{itemize}
Moreover, $h^{-1}$ extends continuously to a map $g: \hat{\strip} - \slits \>
\closure{S}$, and there exists $\tilde g$ such that 
$\tilde{g}: \hat{\strip} \> \closure{S}$  is a homeomorphism
with $\tilde g \on
\strip$ quasiconformal,
and $\tilde g \on {\partial \strip} =
g \on {\partial \strip}$
\end{lem}
We will also show:

\begin{lem} \label{slitbounds}
There exists a $B$ such that for all continuous $f: \hat{\strip}
- \slits \> \reals$ with $\norm{f \on {\strip -
\closure{\slits}}}_{1,2} \le 1$,  there exists $\tilde{f}$ on $\hat{\strip}$ with
$\tilde{f}=f$ on $\partial \strip$, $\tilde{f}$ harmonic on
$\strip$, and $\norm{\tilde{f}} \le B$.
\end{lem}

\begin{remark}
The statement of this lemma is a little peculiar, because $f$ is
assumed continuous on a set that is neither open nor closed. It is
certainly not enough to assume that $f$ is continuous on 
$\hat{\strip}
\setminus
\closure{\slits}$. 
\end{remark}

Given these two lemmas, we can quickly prove:
\begin{lem}
$$
\distS {\closure N} S < \infty.
$$
\end{lem}
\begin{proof}
If $f$ is continuous on $\closure S$ with $\norm{f \on {S-N}} \le 1$,
	then $f \circ g$ (with $g$ as in Lemma \ref{mapping})
	satisfies the hypotheses of Lemma \ref{slitbounds},
	so we can find $\widetilde{f \circ g}$ with $\norm{\widetilde{f \circ g} \on \strip}
\le B$,
 and then $\tilde f:= \widetilde {f \circ g} \circ \tilde{g}^{-1}$ has universally
bounded Sobolev norm on $S$
 (by Fact \ref{q q}), and $\tilde f \on {\partial S} = f \on \partial S$.
\end{proof}

\subsection{Proof of Sobolev bounds for the slitted strip ($\strip, \slits$).}

\setproof{ of Lemma \ref{slitbounds}}{\ref{slitbounds}}
\begin{proof}
We first need to describe a formula for the $\Sobolev$ norm of a
harmonic function on $\strip$. Let $\blackboard H$ denote the upper
half plane. Suppose that $g:\reals \cup \{\infty\} \>\reals$ is
continuous, and continuous at infinity, in the sense that $\lim_{t\>
\infty} g(t)$ exists (and is independent of direction). Then there is
a unique continuous harmonic extension $\tilde{g}$ of $g$ to
$\blackboard H$, and its Sobolev norm on $\blackboard H$ is given by
$$
\norm{\tilde g}^2 = \frac 1{2\pi}\int_{-\infty}^{\infty}\!
\int_{-\infty}^{\infty}
\frac{(g(s) -g(t))^2}{(s-t)^2} ds \, dt.
$$
 (In particular, $\tilde g \in \Sobolev$
if and only if the double integral is finite.)
This formula appears as equation (24) in \cite{Nag:Sullivan:quantum}.
\crit{Say also, ``Look also at equation (20) there for their 
normalization of the Sobolev norm'', or normalize the norm the way that 
they do, and fix the formula}

Now let $f$ be a (real-valued) function on the ideal boundary
$\partial \strip$ of the infinite strip $\strip$; for $t \in \reals$ we let
$f_0(t)=f(t)$, and $f_1(t)=f(t+i\pi$. We require that $f$ is
continuous; this is the same as saying that the $f_i$ are continuous,
and that $\lim_{t \> +\infty} f_i(t)$ and $\lim_{t \> -\infty} f_i(t)$ each
exist and are independent of $i$.  Using the conformal map $z \mapsto
e^z$ from $\strip$ to $\blackboard H$, we obtain the following
formula:
\bl Let $f$ on $\partial \strip$ be continuous; then $f$ has a unique
continuous harmonic extension
$\tilde{f}$ to $\strip$, whose Sobolev norm is given by
$$
\norm {\tilde{f}}^2 = \sum_{i,j=0,1} \frac {I_{ij}(f)}{2\pi}, 
$$ where\mpar{I hope this is correct}
$$ I_{ij}(f) =  \int_{-\infty}^{\infty}\! \int_{-\infty}^{\infty}
\frac{(f_i(s) -f_j(t))^2}
{(e^\frac{s-t}{2}-(-1)^{i+j}e^{-\frac{s-t}{2}})^2}ds \, dt
$$
\el
Note that each $I_{ij}(f)$ above is non-negative, so each $I_{ij}$ must
satisfy $I_{ij}(f) \le
2 \pi \norm {\tilde{f}}^2$

Suppose, as in Lemma \ref{slitbounds}, that $f$ is defined and continuous on
$\hat{\strip}
\setminus
\slits$, and $f$ has Sobolev norm at most 1 on $\strip \setminus
\closure{\slits}$. The conditions\newcrit{which conditions?} on $\slits$
imply\newcrit{how?} that
$$ f(t+i\pi)-f(t-i\pi) = \int_{-\pi}^{\pi} \frac{\partial
f(t+iv)}{\partial v} dv \qquad(*)
$$ for almost every $t$. By the Lemma above, to find a bound for the
Sobolev norm of the harmonic extension of $f \on {\partial \strip}$,
we just need to establish bounds for each $I_{ij}(f)$.

\def\nftwo{\norm{f}{}^2}
Let $v: \strip \> \strip$ be defined by $v(x+iy) = x +iy/5$. Then $v$
is 5-qc, and $v(\strip)$ is a substrip $\cal E$ of $\strip$ that lies below
$\slits$ (${\cal E} = \{ z |0 < \Im z < \pi/5\}$). Then $f \circ v \on \reals= f \on
\reals$, so $I_{00}(f \circ v) = I_{00}(f)$.  By Fact \ref{q q},
$\norm{f \circ v}_\strip \le 5 \norm{f}_{\cal E}$. So we obtain:
$$
\frac1{2\pi} I_{00}(f) = \frac1{2\pi} I_{00}(f \circ v) \le
\norm{f \circ v}_\strip^2 \le 25 \norm{f}_{\cal E} \le \nftwo,
$$ and likewise for $I_{11}$.

So we just need to bound $I_{01}$. From the inequalities
$$ 
(f_0(s)-f_1(t))^2 \le 2((f_0(s)-f_0(t))^2 + (f_0(t)-f_1(t))^2)
$$ 
and
$$
\frac1{(e^\frac{s-t}{2}-e^{-\frac{s-t}{2}})^2} \ge
\frac1{(e^\frac{s-t}{2}+e^{-\frac{s-t}{2}})^2}
,
$$ 
we obtain
$$ 
I_{01} \le 2I_{00} + \int_{-\infty}^{\infty}\frac{1} {(e^\frac
s2+e^{-\frac s 2})^2}ds
\int_{-\infty}^{\infty}(f_0(t)-f_1(t))^2 dt.
$$  
Now
\begin{eqnarray*}
\int_{-\infty}^{\infty}(f_0(t)-f_1(t))^2 dt &=&
\int_{-\infty}^{\infty}\left(\int_0^\pi \frac{\partial f(t+iv)}{\partial v}
dv\right)^2 dt \\ 
& &\text{(by $(*)$)} \\
&\le & \pi\int_{-\infty}^{\infty}\left(\int_0^\pi
\left(\frac{\partial f(t+iv)}{\partial v}\right)^2 dv\right) dt \\ 
& &\text{(by the Cauchy-Schwarz inequality)}\\
& \le & \pi\nftwo.\\
\end{eqnarray*}

Thus each $I_{ij}$ is bounded in terms of $\nftwo$, 
so we have bounded $\norm{\tilde f}$.
\end{proof}

\subsection{Proof of mapping lemma}


In order to prove Lemma \ref{mapping}, we first introduce another
canonical object, the ``recursively slitted square''. We show that
there is a map from the recursively notched square to the recursively
slitted square with properties analogous to that described in Lemma
\ref{mapping}.  Then we describe a quasiconformal map from the
recursively slitted square to the strip $\strip$ that maps the slits
of the recursively slitted square to a union $\slits \subset \strip$ of
slits with the properties described in Lemma \ref{mapping}.


\def\sqslits{\slits'}
Let us now define the recursively slitted square. 
Let $S'$ denote the open square $(-1,1) \times (-1,1)$.  We now
define a set $\sqslits \subset S'$, which is the union of a set of
vertical slits.  Let ${\bf Q\/}_2$ denote the set of dyadic rational points in
the interval $(-1,1)$.  For each $\alpha \in {\bf
Q\/}_2$, let $v_{\alpha}$ be the minimal power of $k$ such that
$\alpha=p/2^k$.  Let $V_{\alpha}$ be the vertical segment given by
$$x=\alpha; \hskip 15 pt |y| \leq \frac 35 2^{-v_{\alpha}}.$$
Define
$$\sqslits=\bigcup_{\alpha \in {\bf Q\/}_2} V_{\alpha}.$$

For future reference (in the proof of Lemma \ref{diamond})
, we note the following:
\begin{fact} \label{angle}
If $x+iy \in \slits'$, then  $|y/(1+x)| \le \frac35$ and $|y/(1-x)| \le \frac35$.
\end{fact}
\begin{proof}
We have $x = p2^{-k}$ with $-2^k < p < 2^k$, $p \in \integers$ , and
$|y| < \frac 35 2^{-k}$. Therefore $x \ge (1-2^k) 2^{-k}$, so $1+x
\ge 2^{-k}$, so $|y/(1+x)| \le \frac35$. Likewise $|y/(1-x)| \le \frac35$.
\end{proof}

\begin{prop} \label{square smashing}
There is a continuous map $\phi: \closure S - \Int N\> \closure {S'}$
with the following properties:
\begin{enumerate}
\item
$\phi(S-\closure N)= S' - \closure {\slits'},$ and $\phi:S-\closure N \> S' -
\closure {\slits'}$ is a quasiconformal homeomorphism.
\item
$\phi(S-N) = S' - \slits',$. and $\phi: S-N \> S' - \slits'$ is a homeomorphism. 
In particular, $(\phi\on{S-N})^{-1}: S'-\slits' \> S-N$ is continuous.
\item
There is homeomorphism $\psi:\closure S \> \closure {S'}$ such that
$\psi: S \> S'$ is quasiconformal, and $\psi \on {\partial S} = \phi
\on {\partial S}$.
\end{enumerate}
\end{prop}

From what one can tell from word of mouth and Yoccoz's lectures, 
a similar proposition is used by Yoccoz
\cite{Yoccoz:local:connectivity} in his proof of Theorem
\ref{mlc} (local connectivity of $M$ at $c$). 
Yoccoz uses it to prove a more limited version of the quasiconformal
distortion bounds for the recursively notched rectangle.  It seems to
be a folk result: the author is unsure of its original discoverer. The
idea of it was described to him by his advisor, Curtis McMullen. Since
it does not appear in the present literature, we will give a complete
proof of it.

\begin{proof}

The idea of the proof is to divide $S-\closure N$ and $S'
-\closure {\slits'}$ into a countable collection of similar (in the Euclidean
geometry sense) regions organized in a tree-like fashion, and define a
piecewise linear map from each region in $S-\closure N$ to the
corresponding region in $S' -\closure {\slits'}$.  We then check that these
piecewise linear maps fit together to a quasiconformal map
$\phi:S-\closure N \> S' -\closure {\slits'}$, and that $\phi$ extends
continuously to $\closure S - \Int N$, and that the extension has the
desired properties.

We say that a {\it marked rectangle\/} is a pair
$(A,B)$ where $B$ is a rectangle and 
$A \subset \partial B$ is closed subinterval properly
contained in a side of $\partial B$.  We say that
a pair $(A',B')$ is a {\it slitted rectangle\/} if
$B'$ is a rectangle, and $A' \subset B'$ is a
segment perpendicular to $\partial B'$ which 
intersects $\partial B'$ in a single point. (See figure \ref{0}).

\begin{figure}
\centerline{\epsfbox{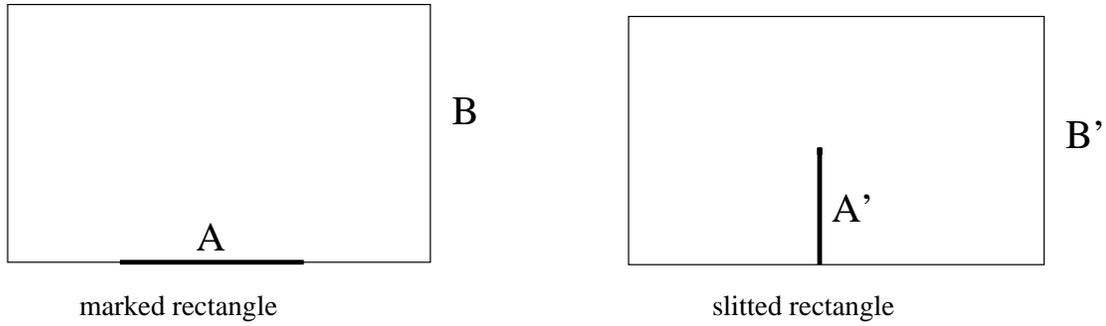}}
\caption{The marked rectangle and slitted rectangle.}
\label{0}
\end{figure}

The two combinatorially equivalent triangulations shown in figure
\ref{1} determine a piecewise affine (and hence quasiconformal) map
$\alpha: B-A \to B'-A'$, defined by letting $\alpha$ on each triangle
be the unique affine map mapping the triangle to the corresponding
primed triangle.  This PL map $\alpha$ will be the building block for
the desired quasi-conformal map from $S-\closure N$ to $S' -
\closure{\slits'}$.

\begin{figure}
\centerline{\epsfbox{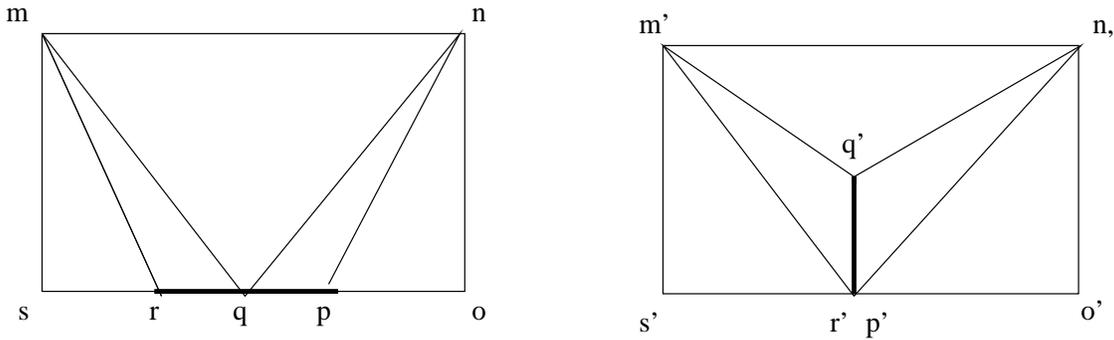}}
\caption{Combinatorially equivalent triangulation of the 
marked and slitted rectangles.}
\label{1}
\end{figure}

Let $X$ denote the union of horizontal lines
in the plane of the form $y=\pm \frac 12 3^{-n}, n >0$.
  Figure \ref{2} shows how
the lines of $X$ intersect $S - \closure N$ and partition it into
connected components.  Each component of the partition is a marked
rectangle, and the components are in fact all similar to each other.

\begin{figure}
\centerline{\epsfbox{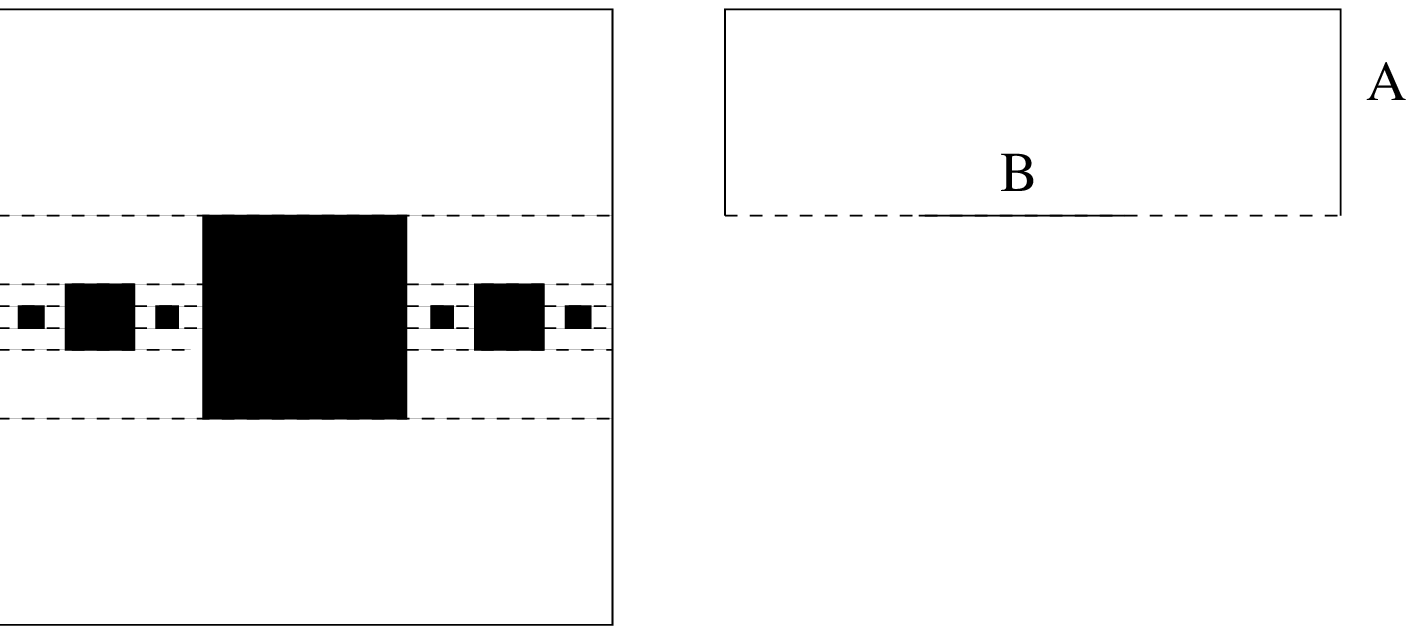}}
\caption{How the lines of $X$ intersect $S - \closure N$.}
\label{2}
\end{figure}

Let $X'$ denote the union of horizontal lines of the
form $y=\pm 2^{-n}$.  Figure \ref{3} shows
how the lines of $X'$ intersect 
$S -\closure{\slits'}$ and partition it into connected
components.  Each component is a slitted rectangle,
and the components are all similar to each other.

\begin{figure}
\centerline{\epsfbox{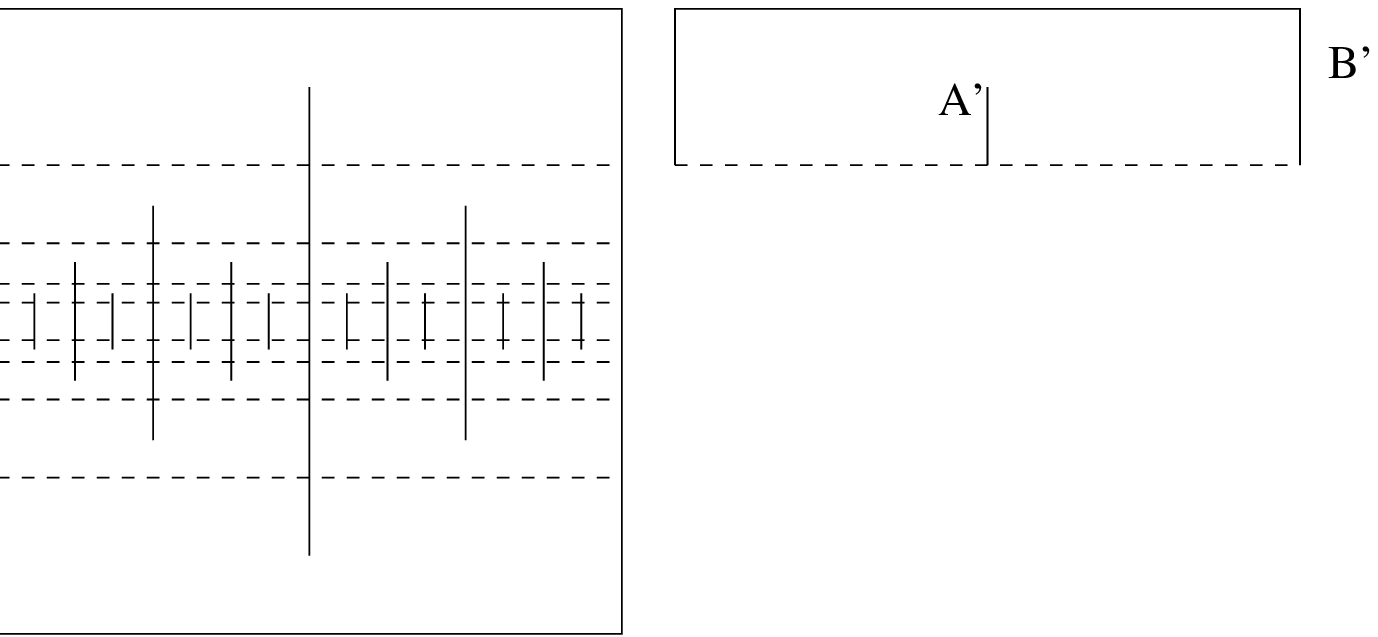}}
\caption{How the lines of $X'$ intersect $S' -\closure{\slits'}$.}
\label{3}
\end{figure}

The components of the partition of $S - (\closure N \cup X)$
correspond bijectively with the nodes of a pair of infinite binary
trees---one for the top half of $S - \closure N$ and one for the
bottom.  Likewise for the components in the partition of $S
-\closure{\slits'}$.  In fact, the combinatorial structure of the two
partitions is the same.  The map $\alpha$ defined above extends,
component by component, to give a piecewise linear, quasi-conformal,
map $\phi: S-\closure{N} \to S' -\closure{\slits'}$.  We will verify, in
turn, that $\phi$ has the properties stated in Lemma \ref{square
smashing}.

Observe that $\phi$ extends continuously to the union of the closures
of the marked rectangle components of $S-\closure N - X$. This union
is equal to $\closure S-(\Int N \cup \reals)$.
So we just need to check that it extends continuously to $\reals - \Int N$, 
which is just the middle-thirds Cantor set $\mt$ \newcrit{$\mt \neq M$, the
Mandelbrot set}.

\crit{Rewrite this to make it more direct---we've already defined $M$}
 The ends \mpar{ends of the other binary tree}
of the pair of binary trees for the
partition $S-(\closure N \cap X)$ correspond to the middle thirds
Cantor set in $S$, consisting of
all points of the form $\sum_{i=1}^\infty a_i 3^{-i}$ with each
$a_i\in \{0,2\}$. Therefore\newcrit{?} $\phi$, so far defined on $S- \Int N -\reals$,
extends continuously to this Cantor set subset of the reals as the
Cantor function $\sum_{i=1}^\infty a_i 3^{-i} \mapsto 
\sum_{i=1}^\infty \frac{a_i}2 2^{-i}$.

 Thus we have defined a continuous map
$\phi: S-\Int N \> S'$, and we have already seen that property 1, that
$\phi:S-\closure N \> S'-\closure \slits$ is a quasiconformal
homeomorphism, is satisfied.

Property 2 then follows from the following simple lemma in point-set topology:

\begin{lem}
Suppose there exist $X$,$Y$,$f:X \> Y$,and $A \subset X$ such that
\begin{enumerate}

\setlength{\parskip}{-3pt plus 3pt minus 1pt}

\item
$X, Y$ are compact metric  spaces,
\item
$f:X \> Y$ is continuous,
\item 
$f \on A$ is injective, and $f(A) \cap f(X-A)= \nullset$.

\end{enumerate}
Then $f\on A: A \> f(A)$ is a homeomorphism.
\end{lem}

\goodbreak

Note that we do not assume that $A$ is a closed subset of $X$. We
could drop the condition of metrizability, at the expense of using
nets in the proof instead of sequences. 

\begin{proof}
Note that $f^{-1}$ is a well-defined function on $f(A)$. We just need to
show that it is continuous, which is equivalent to showing that if
$y_i, y \in f(A)$, with $\lim_{i\>\infty} y_i \> y$, then $f^{-1}(y_i)
\> f^{-1}(y)$. It is enough to show that every subsequence of the
$y_i$ has a subsequence with the above property (that $f^{-1}(y_i) \>
f^{-1}(y)$---here we follow the convention of not changing notation
for passing to subsequences). So, given a subsequence of the $y_i$,
pass to a further subsequence such that $f^{-1}(y_i) \> z$ for some $z
\in X$ (possible by the compactness of $X$). 
But then $f(z) = y$ by the continuity of $f$, which implies that $z$
is equal to $f^{-1}(y)$, the unique element of $f^{-1}(\{y\})$.
\end{proof}

So we just apply this Lemma to the case where $X=\closure S - \Int N$, $ 
Y=\closure{S'}$, $ f = \phi:\closure S - \Int N \>\closure{S'}$,
and $ A=\closure S-N$, and thereby conclude that 
$\phi:S -N \> S'-\slits'$ is a homeomorphism.

\medskip

Finally, to show property 3, consider the combinatorially equivalent
partitions of $S$ and $S'$ depicted in Figure \ref{4}. The crescent
shaped sets of each partition (there are countably many---only
finitely many are shown in the figure, of course) are equivalent to
each other by Euclidean similarities.  In $S$, their sizes decrease in
powers of $3$.  In $S'$, their sizes decrease in powers of $2$.  Each
piece in the partition of $S$ may be mapped to the piece in the
partition of $S'$, in the same combinatorial location, by a piecewise
linear map, whose dilatation is independent of the choice of piece in
the partition.  The resulting map $\psi$ is therefore quasi-conformal.
It agrees with $\phi$ on $\partial R$, because the boundary values of
both maps are piecewise linear maps (on $\partial S-\closure N$---note
here that $\partial S
\cap \closure N$ consists of just two points) that respect the same partitions
of $\partial S-\closure N$. 
\crit{Fix up language here a bit}
\begin{figure}
\centerline{\epsfbox{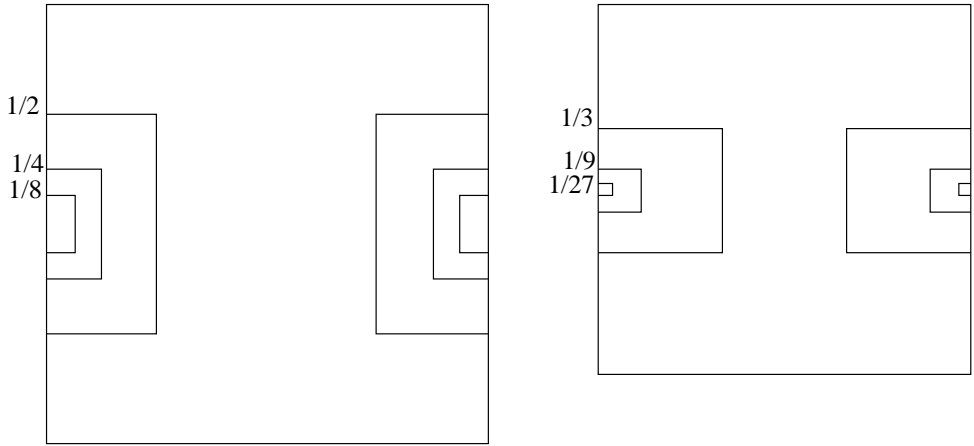}}
\caption{Combinatorially equivalent partitions of $S$ and $S'$.}
\label{4}
\end{figure}

This completes the proof of the proposition.
\end{proof}


We now describe a quasiconformal map from the recursively slitted square to a
``ruler'' (defined below), thus completing the proof of Lemma \ref{mapping}.

Let\mpar{Rich's contribution}
 $\strip$ be the infinite strip given by $|\Im z| \leq 1$.
Let $\slits \subset \strip$ be the union of a countable collection of
vertical intervals, having the following properties:
\begin{enumerate}
\item   For all $z \in \slits$, $|\Im z| <\frac 35$
\item The set $\slits$ is symmetric with respect 
to reflection in the $x$-axis. 
\item $\closure \slits \subset \slits \union \reals$
\end{enumerate}
We say that the pair $(\slits, \strip)$ is a {\it ruler\/}.  The value
$\frac 35$ above is taken to correspond to the requirement that all
the slits composing $\slits$ in the statement of Lemma \ref{mapping} have
imaginary part between $\pi/5$ and $4\pi/5$ (in the strip defined by $0
< \Im z < \pi$.)\newcrit{Should say how to get from one coordinate system to another}

\begin{lem}\label{diamond}
There exists a quasi-conformal map $\phi: S' \to \strip$ with
the following properties:
\begin{enumerate}
\item $\phi$ is symmetric with respect to reflection in the
coordinate axes. 
\item $\phi$ takes $(\slits',S')$ to a ruler
$(\slits,\strip)$.
\end{enumerate}
\end{lem}

\begin{proof}
Let $Q$ be the ``diamond'' inscribed in $S'$, whose vertices
are at the midpoints of sides of $S'$.  There 
is a simple quasi-conformal (in fact, piecewise-linear) map from $S'$ to $Q$ 
(See figure \ref{clear}---we map $\triangle acd$ to $\triangle ac'd$ and
$\triangle bcd$ to $\triangle bc'd$, and likewise in the other three corners. 
The map is in fact the identity on the convex hull of $\slits'$. Fact \ref{angle}
implies that the convex hull of $\slits'$ is indeed as shown in figure \ref{clear}) 
 Thus,
we may work with the pair $(\slits',Q)$ instead of with
$(\slits',S')$.  
Let $Q_-$ denote those points of $Q$ with negative
$x$ coordinate.  Likewise define $Q_+$. 

\begin{figure}[htb]
\centerline{\epsfbox{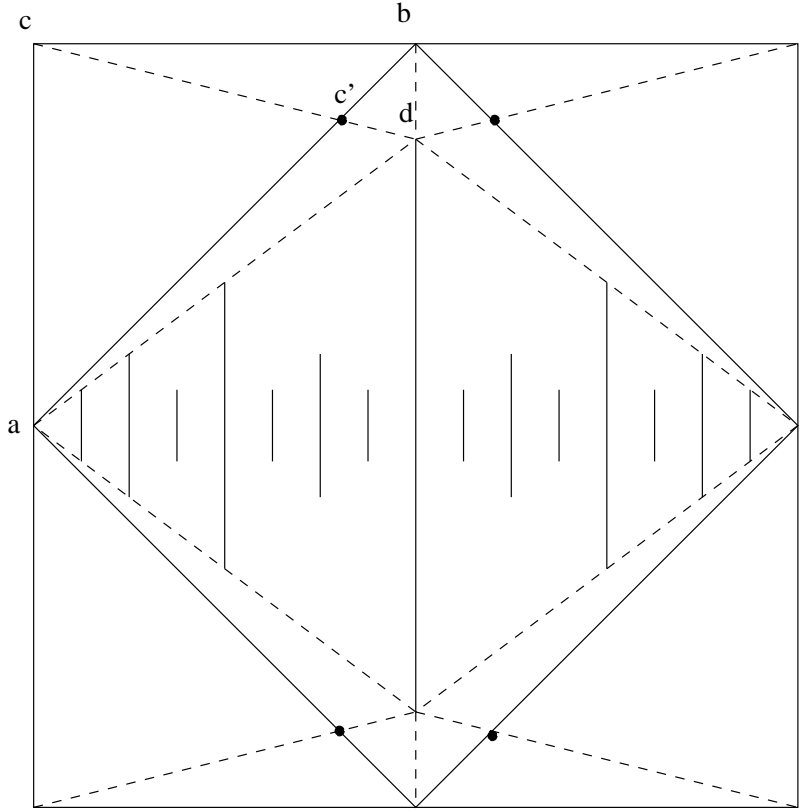}}
\caption{Piecewise linear map from $S'$ to $Q$ fixing $\slits'$ pointwise}
\label{clear}
\end{figure}

Define $\rho_-: Q_- \to R'$ by the formula
$$\rho_-(x,y)=(\log(1+x),y/(1+x)).$$
Define $\rho_+: Q_+ \to R'$ by the formula
$$\rho_+(x,y)=(-\log(1-x),y/(1-x)).$$
We compute the Jacobian:
$$D\rho_-=
\left[
\matrix{
1/(1+x) & -y/(1+x)^2 \cr
0 & 1/(1+x)}
\right].$$
Multiplying by $1+x$, we have:
$$(1+x)D
\rho_-=
\left[
\matrix{
1 & -y/(1+x) \cr
0 & 1}
\right].$$

Finally, by Fact \ref{angle}, $|y/(1+x)| \leq 1$ in $Q_-$. Therefore (1+x)$D
\rho_-$, stays within a compact subset of ${\blackboard{GL_2(R)}}$, so
the dilatation of $D \rho_-$, which is equal to the dilatation of
$(1+x)D\rho_-$, is uniformly bounded.  (In fact, $\norm{(1+x)D
\rho_-} \le \sqrt 3$ (because the square of the norm of a matrix is less
than the 
sum of the squares of its entries), and $(1+x)D\rho_-$ is area-preserving, 
so the dilatation of $(1+x)D \rho_-$ is at most 3). A similar computation
proves the same for $\rho_+$.  Note that both these maps take vertical
line segments to vertical line segments.

The union
$\rho=\rho_- \cup \rho_+$ defined on $Q_- \cup Q_+$ is symmetric with respect to
reflection in the $y$-axis, and hence extends to
a quasiconformal map on all of $Q$.

Finally, we must check that all points in $\rho(\slits')$ have absolute
value of imaginary part less than $\frac 35$. It is enough to check this for 
$\rho_-(\slits' \cap Q_-)$. Suppose $x+iy \in \slits' \cap Q_-$. Then 
$| \Im \rho_-(x+iy)|=
|y/(1+x)| \le \frac 35$ 
by Fact \ref{angle}.  
\end{proof}

So, to prove Lemma \ref{mapping}, simply follow the map given in
Proposition
\ref{square smashing} 
with the map given by Lemma \ref{diamond}, and then map the 
resulting strip by a Euclidean similarity to the one described for
Lemma \ref{mapping}.
\crit{Say a bit more here?}

\def\ray#1{R(#1)}
\def\landing#1{l(#1)}
\def\rayland#1{\landing{\ray {#1}}}
\def\grp#1#2{\closure R(#1,#2)}
\let\rp=\grp
\def\crp#1#2{\meet{#1}{#2}, #1 \neq #2}
\def\eq{\simeq}
\def\meet#1#2{#1 \simeq #2}
\def\t{\theta}
\def\cara{Carath\'eodory}
\section{Covering $J$ with the image of the recursively notched square.}
\label{embedding rns}

\newcrit{Talk explicitly about equivalence and {\em emphasis\/}[?]:  rays meeting
in pairs/the equivalence induced by $\phi \on {S^1}$}

The purpose of this section is to prove Lemma \ref{getS}, which says
roughly that we can embed copies of the recursively notched square
into each piece $P$ so as to cover the Julia set near the boundary of
the piece with copies of $\closure N$. The embeddings are to be
quasiconformal on $S-\closure{N}$.  We first do so in the case where
$P$ is the top level piece containing the critical value, in which case we need
only embed one copy of the RNR. We then use that embedding and the
dynamics to get embeddings for all other pieces.

To get the embedding for the top level piece containing the critical value,
 we proceed
as follows.  Denote the two arguments of the external rays bounding that
piece by $A$ and $D$, with $A < D$. Then we find intervals $[A,B]$ and
$[C,D]$ and a pair of monotonic maps $q_1:\mt \> [A,B]$ and $q_2: \mt \>
[C,D]$ such that, for each $x \in \mt$, the rays with arguments $q_1(x)$
and $q_2(x)$ land at the same point in $J$. \newcrit{Recall that $M \in [0,1]$ is
the middle-thirds Cantor set.} Thus $\phi \circ e^{2 \pi i
q_i(\cdot)}: \mt \> J$ is independent of $i$. \newcrit{$\phi$ is defined below
(or in section 1.3)---maybe we could consolidate into a single background section}
 The maps $q_1,q_2$ will also be such that $\phi \circ e^{2 \pi i q_i}$ will
extend to a quasiconformal map of $\closure S - \Int N$ into the
dynamical plane. \newcrit{Actually a homeo of $\closure{S}-N$ that is qc on
$S-\closure N$.} That map can then be easily extended to a map of
all of $S$ into the dynamical plane, such that $N$ covers a neighborhood of
$\partial P$ in $J$. 
\crit{The point really is not that $N$ is covering, 
but that $\closure S-\closure N$ avoids the Julia set.}

\newcrit{Maybe put in some sort of subsection-by-subsection plan/outline.}

\crit{Maybe also give a section-by-section overview}
\subsection{Definitions and observations for external rays}
We first require some basic definitions. \newcrit{Also recall $\alpha, \beta$.}

Recall\mpar{Some of this is in 1.3, some of it isn't}
that there exists a unique conformal isomorphism $\phi: \cx -
\closure{\Delta} \> \cx - J$
such that, for all $z \in \cx -\Delta$,
$\phi(z^2)= (\phi(z))^2+c$.
Because $J$ is locally connected (Theorem \ref{jlc}), \cara 's Theorem
implies that $\phi$ extends continuously to a map
$$
\phi: \cx - \Delta \> \cx
$$
so that $\phi(\partial \Delta) = J$. \newcrit{The extension also satisfies the
functional equation}

Recall also that a \term{external ray} (or just \term{ray}) $\ray \theta$ is
defined by 
$$
\ray \theta := \{\phi(r  e^{2 \pi i \theta})  \mid 1< r  <\infty \}.
$$
Here we think of $\t$ as an element of $\reals/\integers$. Each such element has
a unique representative in $[0,1)$; we may sometimes denote the element by such a
representative. The conjugacy properties of $\phi$ imply that $f(\ray \theta)=
\ray {2\theta}$. 

We say that a ray \term{lands} at $z \in J$ if 
$$
\lim_{r \> 1} \phi(r  e^{2 \pi i \theta}) = z,
$$
which is equivalent to saying that 
$\closure {\ray \theta} = \ray \theta \cup \{z\}$.
Because $J$ is locally connected, 
\cara 's theorem 
\crit{state it?}
 implies 
that every ray $\ray \theta$ lands.
We denote the landing point of $\ray \theta$ by $\landing {\ray \theta}$.
\cara 's theorem also implies that $\landing {\ray \theta}$ is a 
continuous function of $\theta$.  Note that $\rayland \t =
\phi(e^{2\pi i\theta})$. \newcrit{Some redundancy in notation here.}
If $\rayland {\t_1} = \rayland{\t_2}$, then
we write $\t_1 \eq \t_2$.


The term \term{combinatorial ray-pair} will denote a pair $(\theta_1,
\theta_2)$ such that 
$\theta_1 \neq \theta_2$ but
$\t_1 \eq \t_2$. 
The term \term{geometric ray-pair} will denote the union $\ray
{\theta_1} \cup \ray {\theta_2} \cup \{z\}$, where $(\theta_1,
\theta_2)$ is a combinatorial ray-pair, and $z=\rayland {\theta_1} =
\rayland{\theta_2}$.  
We will denote the geometric ray-pair corresponding to the
combinatorial ray-pair $(\theta_1, \theta_2)$ by $\rp
{\theta_1}{\theta_2}$.
We will use the term \term{ray-pair} to refer to either a
combinatorial or geometric ray-pair when the context makes it clear
which one is being referred to.

Note that the continuity of $f$ and conjugacy properties of $\phi$
imply that if $\rayland \theta = z$, then $\rayland {2\theta} = f(z)$. 
Therefore, if $\t_1 \eq \t_2$ then 
$2 \t_1 \eq 2 \t_2$.
The converse, however, is not true:
$\meet{2 \theta_1}{2 \theta_2}$ 
does not necessarily imply $\meet{\theta_1}{\theta_2}$.
We need to describe circumstances in which some sort of partial converse 
can be obtained. \crit{Clarify?}

\begin{defn}
A \term{slice} is an open subset $S$ of $\cx$ such that the boundary
of $S$ has two components, each of which is a geometric ray-pair.
\end{defn}

Any two distinct geometric ray-pairs bound a unique slice.  The slices
$S$ we will be interested in have the property that  $\rayland 0
\notin S$.  \newcrit{Recall that $\beta=\rayland 0$.}
Such slices are called
\term{vertical slices}. \newcrit{Maybe explain why they're called
\term{vertical}}  We can write the boundary of a vertical slice
$S$ as
$\grp ad \cup \grp bc$, where $0 < a < b < c < d < 1$.  In this case we write
$S=S(a,b,c,d)$.

\def\cong{\equiv}
\begin{lem} \label{pair pullback1}
If $f^n: S(a,b,c,d) \> S(a',b',c',d')$ is univalent, and
$2^n a \cong a',2^n b \cong b', 2^n c \cong c',2^n d \cong d'$ (all modulo 1), 
then $b'-a'= 2^n(b-a), d'-c' = 2^n(d-c)$\newcrit{literally, not modulo 1}, and, if
$x,y$ are such that
$a < x < b$ and $c < x < d$, and $2^n x \eq 2^n y$, then $x \eq y$.
\end{lem}

\begin{proof}
\crit{Show $b'-a'= 2^n(b-a), d'-c' = 2^n(d-c)$.}
We know that $f^n$ is injective on the rays in $\closure{S(a,b,c,d)}$, so
$z \mapsto 2^n z$ is injective on $[a,b] \cup [c,d]$, so we have
$b'-a'= 2^n(b-a), d'-c' = 2^n(d-c)$.

We have $a' < 2^n x < b'$ and $c' < 2^n y <
d'$. Therefore, $\grp {2^nx}{2^ny} \subset S(a',b',c',d')$. Denoting
the inverse of $f^n: S(a,b,c,d) \> S(a',b',c',d')$ by $g$, we find
that $g(\grp {2^nx}{2^ny})$ must be a geometric ray-pair $\grp uv$
with $2^nu \cong 2^nx,2^nv \cong 2^ny$, and $a < u < b, c < v < d$. Since $b-a <
2^{-n}$ and $d-c < 2^{-n}$, we must have $u=x$ and $v=y$. Therefore $x
\eq y$.
\end{proof}

We can also state the analogous result 
when the slice is ``flipped over'' by $f^q$.
(that is, when the ``vertical orientation'' is reversed). 
\begin{lem} \label{pair pullback2}
If $f^n: S(a,b,c,d) \> S(a',b',c',d')$ is univalent, and
$2^n a=c',2^n b=d', 2^n c=a',2^n d=b'$ (all modulo 1), 
then $b'-a'= 2^n(d-c), d'-c' = 2^n(b-a)$, and, if $x,y$ are such that
$a < x < b$ and $c < x < d$, and $2^n x \eq 2^n y$, then $x \eq y$.
\end{lem}

The proof is the same.

\subsection{Getting univalent slice dynamics}
\label{get slices}

\newcrit{Having set up basic notation and lemmas \ref{pair pullback1} and
\ref{pair pullback2}, we now find the slice dynamics in the top-level piece
containing $c$.}

\begin{figure}
\centerline{\epsfbox{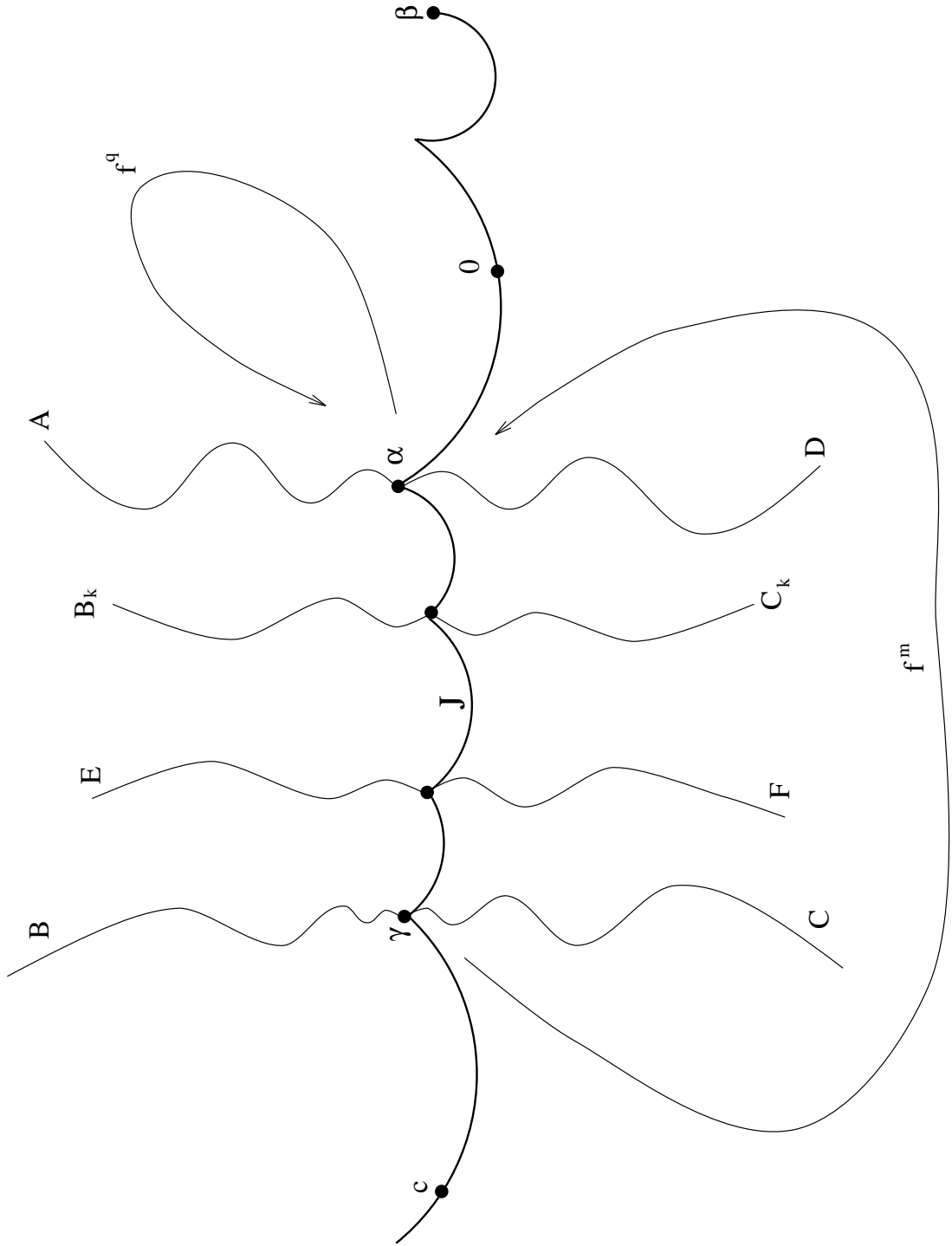}}
\caption{Slices in the dynamical plane}
\label{slices}
\end{figure}
\def\c{{\cal C}}
In this subsection, the information shown in figure
\ref{slices} is built up.  
The reader may wish to check with that figure while reading what follows. 

Now let $\alpha$ denote the $\alpha$ fixed point, at which $q>1$ rays
land.  Each ray is mapped to itself by $f^q$. These $q$ rays union
$\{\alpha\}$ divide $ \cx$ into $q$ components; the boundary of each
component is a single ray-pair containing $\alpha$.  One of these
components (call it $\c$) contains the critical value $c=f(0)$.  The
$\beta$ fixed point $\beta=\rayland 0$ will not be in this component,
because $\beta$ is in the component containing the critical
point. \crit{More explanation?}  Let the ray-pair bounding $\c$ be
$\grp AD$, where $A < D$. (So $\rayland A=\rayland D=\alpha$.)
Because $\rayland 0 \notin \c$, every ray $\ray \t$ landing at a point
in $\c$ will satisfy $A < \t < D$.

For each $n \ge 0$, there is a level $n$ Yoccoz puzzle
piece $P_n$ touching $\alpha$ (i.e.
$\alpha \in \partial P$) and contained in this component $\c$. (So part of the 
boundary of $P$ is $\grp AD$, cut off by an equipotential). If $n$ is
sufficiently large, then $c \notin \closure{P_n}$ (because the diameter of $P_n$
goes to 0 as $n \> \infty$ (Theorem \ref{dzero})). Choose the least such
$n$.
\crit{We will actually have $n =N+1$. Should I say this?}.
The boundary of $P_n$
consists of portions of an equipotential and a finite set of
ray-pairs, cut off at that potential. 
\crit{Explain this description
of a puzzle piece?}
 One such ray-pair must separate $\alpha$ from $c$; we can 
denote it by $\grp BC$, with $A < B < C< D$. So then we have a vertical slice
$S(A,B, C, D)$. Let $\gamma=\rayland B(=\rayland C)$, 
so $\partial (S(A,B,C,D) \cap J)=\{\alpha, \gamma\}$.

We now find smaller slices and univalent maps with which to apply Lemmas 
\ref{pair pullback1} and \ref{pair pullback2}. 

First note that $\{f^i (0) \mid 0 < i \le q\} \cap S(A,B,C,D) =
\nullset$, because $\{f^i (0) \mid 0 < i \le q\} \cap \c = {c}$, and $c
\notin S(A,B,C,D)$. 
Therefore, we can define a single-valued univalent branch $g$ of
$f^{-q}$ on a neighborhood of $\closure{S(A,B,C,D)}$, with
$g(\alpha)=\alpha$, $g(\ray A)=\ray A$, and $g(\ray D)=\ray D$. Then
define $B',C' \in [0,1)$ such that $g(\ray B)=\ray{B'}$ and $g(\ray
C)=\ray{C'}$. Then $f^q:S(A,B',C',D)
\> S(A,B,C,D)$ is a univalent map of vertical slices, mapping boundary
rays to the corresponding boundary rays, and thus satisfies the
hypothesis of Lemma \ref{pair pullback1}

In fact, if we let $B_k, C_k$ be such that $g^k(\ray B)=\ray {B_k}$, and
$g^k(\ray C)=\ray {C_k}$ (so $B_1=B',\quad C_1=C'$), then we have a series of
vertical slices $S(A,B_k,C_k, D)$, and, for each $k \ge 1$,
$f^{kq}:S(A,B_k,C_k, D)\> S(A,B,C,D)$ is a univalent map satisfying
the hypothesis of Lemma \ref{pair pullback1}.
Moreover, $B_k \> A$ and $C_k \> D$ as $k \> \infty$, so 
the diameters of $S(A, B_k, C_k, D) \cap J$ go to zero as $k \> \infty$. 

Now $f^n(\gamma)=\alpha$, since $\gamma$ belongs to the boundary of a
level $n$ puzzle piece. Furthermore, for some $m \in [n, n+q)$,
$f^m(\ray B)=\ray D$, and $f^m(\ray C)=\ray A$.
\ignore{
Also, $f^m$ is injective on some neighborhood of $\gamma$ in
$\closure(S(A,B,C,D))$. \crit{Say why? Explain the case where
$f^k(\gamma)=0$?}  
} 
If $k$ is large enough, then $\{f^i(0) \mid
0 < i \le m\} \cap S(A, B_k, C_k, D) = \nullset$.  Then we can let $h$
be the branch of $f^{-m}$ defined on $S(A, B_k, C_k, D)$, such that
$h$'s extension (also called $h$) to $\closure{S(A, B_k, C_k, D)}$
satisfies $h(\alpha)=\gamma, h(\ray A)=\ray C, h(\ray D)=\ray B$.
Let $E,F$, $0< E < B < C < F < 1$ be such that $h(\ray{B_k})=\ray
F$, and $h(\ray{C_k})=\ray E$.  Then $f^m: S(E,B,C,F) = S(A,B_k,C_k,
D)$ satisfies the hypotheses of Lemma \ref{pair pullback2}, and so does
$f^{m+kq}: S(E,B,C,F) \> S(A,B,C,D)$.

As $k \> \infty$, $B_k-A$ and $D-C_k$ tend to 0, and therefore $B-E <
 D-C_k$ and $F-C < B_k -A$ also tend to 0. So we can choose $k$ such
 that $B-E+B_k-A < B-A$ and $F-C + D-C_k < D-C$, and thereby obtain
 that $A < B_k < E < B < C < F < C_k < D$.  We have now determined all
 of what is shown in figure \ref{slices}.

\subsection{Mapping the RNS into the slice}
\label{mapping rns}
\newcrit{We're still constructing a conguacy from homothety dynamics (on $\mt$)
to boundary-slice dynamics, and our map is still qc at some places. It's hard to
produce a good title.}
\def\mtop{q_1}
\def\mbot{q_2}
\def\mmtop{q_1}
\def\mmbot{q_2}
\def\qs{q}
\def\qc{Q}
\def\emb{\xi}
\def\through{\psi}

Let $\through: \{z \mid \Im z \ge 0\} \> \cx$ be defined by
$\through(z)=\phi(e^{2 \pi i z})$.
Given $y \in \reals^+$, and a vertical slice S(a,b,c,d) we define the
``cut-off slice'' $CS(a,b,c,d,y)$ by $CS(a,b,c,d,y):=S(a,b,c,d) \cap
\through \{z \mid y > \Im z \ge 0\}$. \newcrit{Thus $CS$ is $S$ cut off by an
equipotential}
 So $CS(a,b,c,d,y)$ will be a
bounded domain, and $J \cap CS(a,b,c,d,y) = J \cap S(a,b,c,d)$.  Also,
$\partial (CS(a,b,c,d,y))$ is piecewise smooth curve, and is thus
holomorphically removable.

Our eventual goal is get a homeomorphism $\emb: \closure S \> \closure
{CS(A,B,C,D,1/2)}$ such that $\emb\on {S -\closure N}$ is quasiconformal,
and $\emb(\closure S - \closure N) \cap J=\nullset$. This will be the embedding
required by Lemma \ref{getS}, in the case of the top-level piece containing the 
critical value. 

\begin{figure}
\centerline{\epsfbox{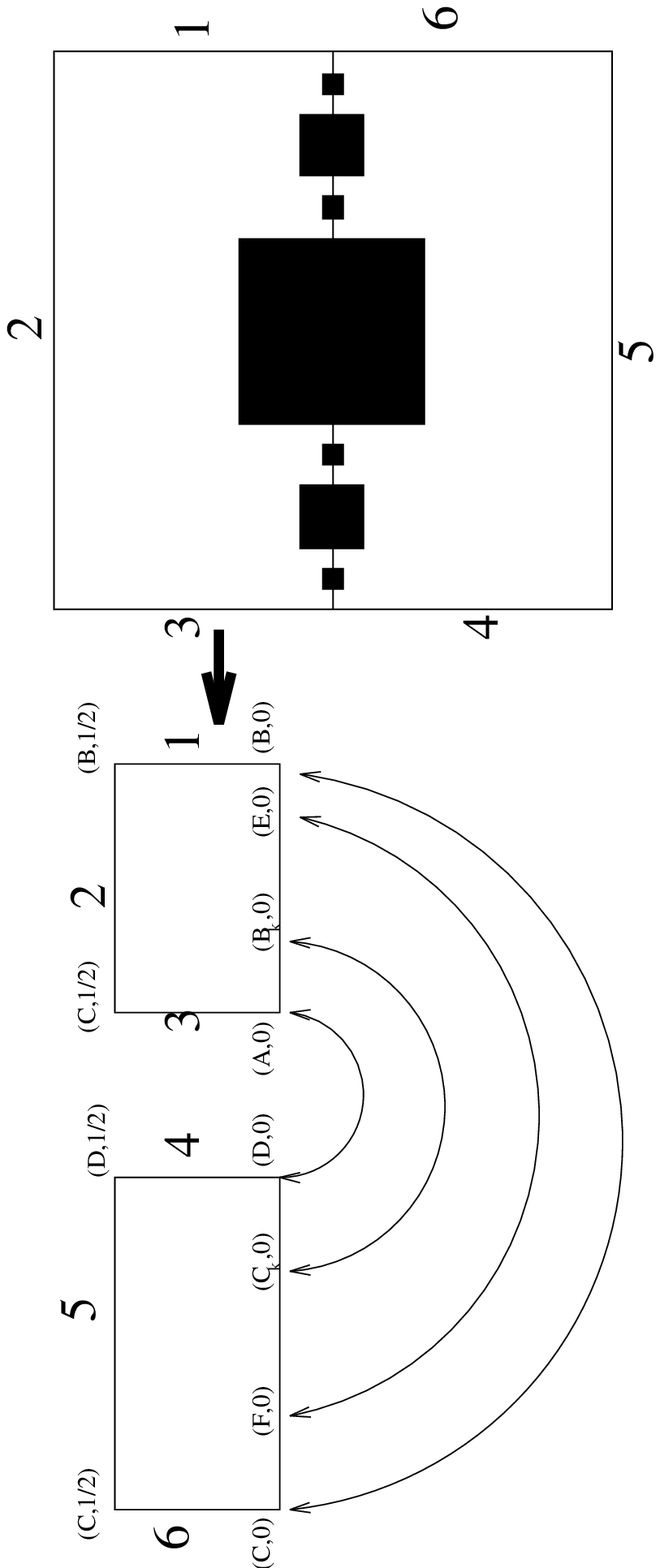}}
\caption{Boundary values for $\mtop$ and $\mbot$
 (numbers 1-6 indicate corresponding sides)
\newcrit{1,2,3 for $q_1$ and 4,5,6 for $q_2$ (so that $q(\partial S)
= \partial \closure{CS(a,b,c,d,1/2)}$).} }
\label{godknows}
\end{figure}

Our plan now is to define embeddings
$$
\mtop: (\closure S -\Int N) \cap \{ z \mid \Im z \ge 0 \} \>
[A,B] \times [0,1/2], 
$$
and
$$
\mbot:  (\closure S -\Int N) \cap \{ z \mid \Im z \le 0 \} \>
[C,D] \times [0,1/2]
$$
that are quasiconformal on the interior of their domains, 
that satisfy 
$\through \circ \mtop \on \mt =\through \circ \mbot \on \mt$
(recall $\mt=(\closure S-\Int N) \cap \reals$), and
that have boundary values as shown in figure \ref{godknows}.
Given such maps, we can then define 
$$
\through \circ (\mtop \cup \mbot): \closure S -\Int N \> \closure{CS(A,B,C,D,1/2)}
\subset  \cx,
$$
which we can then extend to $\closure S$ to get the desired embedding $\emb$. 

Such maps $\mtop$ and $\mbot$ must have the property that for each $x \in \mt$,
$\mtop (x) =(a,0)$ and $\mbot(x)=(b,0)$, with $a \eq b$.
So our goal now is to identify a Cantor set of pairs $(a,b)$. This is done via the
dynamics of slices obtained in the previous section, which we now abstract as
follows: 

On the product of intervals $[A,B] \times
[C,D]$, we then have the following linear dynamics. \newcrit{Confusing---where
does this come from? But comprehensible} Define a linear isomorphism
$$
l_1: [A,B] \times [C,D] \> [A,B_k] \times [C_k, D] \subset  [A,B] \times [C,D]
$$
by the formula
$$
l_1(A+x,D-y)=(A+2^{-kq}x,D-2^{-kq}y).
$$
Then, by Lemma \ref{pair pullback1}, if $(a,b)$ is a combinatorial ray-pair, and
$l_1(a,b)=(a',b')$, then $(a',b')$ is a combinatorial ray-pair, and $g(\grp
ab)=\grp{a'}{b'}$. 
Define another linear isomorphism
$$
l_2: [A,B] \times [C,D] \> [E,B] \times [C,F] \subset  [A,B] \times [C,D]
$$
by 
$$
l_2(x+A,D-y)=(B-2^{-m-kq}y,C+2^{-m-kq}x)
$$
Then, by Lemma \ref{pair pullback2}, if $(a,b)$ is a combinatorial ray-pair, and
$l_2(a,b)=(a',b')$, then $(a',b')$ is a combinatorial ray-pair, and $h(\grp
ab)=\grp{a'}{b'}$. 
Note that for $l_2$, $a'$ depends linearly
on $b$, and $b'$ depends linearly on $a$. So we have two functions
from $[A,B] \times [C,D]$ to itself, with disjoint images. (In fact,
a point in the image of one function cannot share either coordinate with a
point in the image of the other).

Now, note that for any finite sequence 
$i_1, i_2, \dots i_k$, with $i_j \in \{1,2\}$,
we have, by Lemmas 
\ref{pair pullback1} and \ref{pair pullback2} (as observed above), that
$l_{i_1} \circ \dots \circ l_{i_k} (A,D)$ is a combinatorial ray-pair. 
Moreover, if $i_1, i_2, \dots$ is an infinite sequence with $i_j \in \{1,2\}$,
then $\lim_{k \> \infty} l_{i_1} \circ \dots \circ l_{i_k} (A,D)$ exists 
(because the $l_i$ are contracting linear maps) and
is also a combinatorial ray-pair. Let $T \subset [A,B]\times[C,D]$ denote
the set of such pairs.

We now just need to define monotonic embeddings
 $\mmtop: \mt \> [A,B]$ and $\mmbot: \mt \> [C,D]$ such that for all 
$x \in \mt$, $(\mmtop(x), \mmbot(x)) \in T$. (Then we will extend 
$\mmtop$ to a quasi-symmetric map $\mmtop: [0,1] \> [A,B]$, and
then to the desired quasiconformal map 
$\mtop: [0,1] \times [0,1/2] \> [A,B] \times [0,1/2]$ (which is then 
restricted to $(\closure S - \Int N) \cap \{ z \mid \Im z \ge 0 \}$). Likewise
for $\mmbot$).

Now consider\newcrit{is there a better way to say this?} the following artificially
constructed pair of linear isomorphisms, from $[0,1]$ to subsets of
itself. Firstly\newcrit{?}:
$$
e_1: [0,1] \> [0,1/3]
$$
defined by 
$$
e_1(x)= (\frac 13 x),
$$ and secondly:
$$
e_2:  [0,1] > [2/3,1]
$$
defined by
$$
e_2(x) = (1-\frac 13 x).
$$

\ignore{
We claim that the $l$'s and the $e$'s are quasisymmetrically conjugate:
\begin{lem} \label{miracle}
There exists quasi-symmetric maps $q_1:[0,1]\> [A,B]$ and $q_2:[0,1]\> [C,D]$
inducing a homeomorphism 
$q = q_1 \times q_2:[0,1]\times [0,1] \> [A,B] \times [C,D]$
such that for $i \in {1,2}$, $q \circ e_i=l_i \circ q$.
\end{lem}

\crit{Define quasi-symmetric, and figure out how the hell to prove this lemma.}
\crit{Option b: go directly to the piece-wise linear quasiconformal map.
And remember that it doesn't actually have to be defined everywhere.}

It follows that, if we let $\mt \subset [0,1]$ be the middle-thirds Cantor set,
then for all $x \in \mt$, we have that $q(x,x)$ is a combinatorial ray-pair, with
$q$ defined as in Lemma \ref{miracle}. In fact $q(\{(x,x) \mid x \in \mt\})= T$.
}

Then we can define $\mmtop$ and $\mmbot$ by 
$$
(\mmtop (x), \mmbot (x))=\lim_{k \> \infty}l_{i_1} 
\circ \dots \circ l_{i_k} (A,D)
$$
when 
$$
x=\lim_{k \> \infty}e_{i_1} \circ \dots \circ e_{i_k} (0)
$$
where $(i_j)_{j=1}^\infty$ ranges over all possible sequences with 
$i_j \in \{1, 2\}$.

Then $\mmtop:\mt \> [A,B]$ is a monotonic embedding, and furthermore
$\mmtop (\mt) \subset [A,B]$ is a bounded geometry Cantor set,
in the sense that for any sequence $(a_j)_{j=1}^k$, $a_j \in \{0,2\}$,
the ratio 
$$
\mmtop(s)-\mmtop(s+t):\mmtop(s+t)-\mmtop(s+2t): \mmtop(s+2t)-\mmtop(s+3t)
$$
 is bounded,
where $s=\sum_{j=1}^k a_j 3^{-j}$, and $t=3^{-(k+1)}$. This is because the
ratio will always be either $A-B_k:B_k-E:E-B$ or
$C-E:E-C_k:C_k-D$.\newcrit{Justification?} It follows that $\mmtop$ has a 
quasi-symmetric extension
$\mmtop: [0,1] \> [A,B]$ (see the end of \cite{Sullivan:renormalization}
for a discussion\newcrit{Should say just where, if possible}). We can likewise get
a quasi-symmetric extension
$\mmbot: [0,1] \> [C,D]$.

Now we need the following lemma:
\begin{lem} \label{fuck}
If $\qs: [0,1] \> [0,1]$ is a quasisymmetric map, it has a continuous extension
$\qc:[0,1]\times[0,1] \>[0,1]\times[0,1]$ that is quasiconformal
on $(0,1) \times (0,1)$.  (Here we identify $[0,1]$ with $[0,1] \times \{0\}$.) 
We can require $\qc$ to fix each side of the square 
$[0,1]\times[0,1]$ setwise (i.e. map each side of the square to itself). 
\end{lem}

\begin{proof}
We start with $\qc(x,0)=(\qs(x),0)$; we then want to define 
$\qc: \partial( [0,1] \times [0,1]):\> \partial ([0,1] \times [0,1])$ 
so that it is quasi-symmetric (with respect to arc-length); we can then extend
$\qc$ to be quasiconformal on $(0,1) \times (0,1)$ (see 
\cite{Lehto:Virtanen}, section II.6). The right way to define 
$\qc$ on $\partial( [0,1] \times [0,1])$ is by multiple reflection (cf.
\cite{Lehto:Virtanen}, section II.7.2) That is,
we set $\qc (0,x)=(0, \qs (x))$, $\qc(x,1)=(1-\qs(1-x),1)$,
and $\qc (1,x)=(1, 1-\qs(1-x))$. It is then easily checked that
$\qc: \partial( [0,1] \times [0,1]) \> \partial ([0,1] \times [0,1])$ 
is quasi-symmetric. 
\end{proof}

Linearly rescaling the domain and range of \ref{fuck}
and applying it to $\mmtop: [0,1] \> [A,B]$, we obtain
the desired quasi-conformal extension 
$\mtop: [0,1] \times [0,1/2] \> [A,B] \times [0,1/2]$, and likewise
$\mbot: [0,1] \times [-1/2,0] \> [C,D] \times [0,1/2]$.

So then we have
$$
\mtop \cup \mbot: ([0,1] \times [0,1/2]) \cup ([0,1] \times [-1/2,0]) \> 
([A,B] \times [0,1/2])\cup ([C,D] \times [0,1/2]).
$$
From this we can obtain, using the relation
$$
\forall x \in \mt:(\psi \circ \mtop)(x,0)=(\psi \circ \mbot)(x,0),
$$
an embedding
$$
\emb:=\psi \circ (\mtop \cup \mbot):
\closure S - \Int N \> \closure{CS(A,B,C,D,1/2)}
$$
that sends $\partial S$ to $\partial \closure{CS(A,B,C,D,1/2)}$, and
that is quasiconformal on $S - \closure N$.  
(Note that 
$\emb \on {(\closure S - \Int N) \cap \{z \mid \Im z \ge 0\}} = \mtop$, and
$\emb \on {(\closure S - \Int N) \cap \{z \mid \Im z \le 0\}} = \mbot$, and 
$\mtop \on \mt = \mbot \on \mt$, so $\emb$ is well-defined.)
Then $\emb$
applied to the boundary of any component of $N$ is a Jordan curve in $\cx$.  We
can then extend $\emb$ continuously to each component of $\Int N$ via
the Schoenflies theorem, to obtain the desired homeomorphism $\emb:
\closure S \> \closure{CS(A,B,C,D,1/2)}$, with $\emb \on {S-\closure
N}$ quasiconformal, and $\emb(\closure S -\closure N) \cap J =
\nullset$.

\subsection{Embedding RNS's in a arbitrary 
piece $P$ to cover ends of $J \cap P$.}

\def\cstwo{CS_{1/2}}
For economy of space in what follows, let us denote $CS(A,B,C,D,1/2)=\emb(S)$,
by $\cstwo$.
So now we have an embedding $\emb: \closure S \> \cx$ such that 
$\emb (\closure S) = \closure{\cstwo}$, and 
$\emb \on {S -\closure N}$ is quasiconformal, and 
$\emb (\closure S - \closure N) \cap J = \nullset$. We wish to 
show Lemma \ref{getS} for all pieces $P$.

First note that $g^k: \cstwo \> CS(A,B_k,C_k,D, 2^{-(kq+1)})$
is a homeomorphism (review subsection \ref{get slices} for definitions
of $g^k, B_k$, etc.). If we denote the top-level piece containing $c$ by
$P_c(0)$ (see section \ref{notes} for a general discussion of notation), then
we observe that $J \cap (P_c(0) - CS(A,B_k,C_k,D, 2^{-(kq+1)}))$ is 
compactly contained in $P$. Moreover, for $k \ge 1$, and $0 \le t < q$, we have
$(f^t \circ g^k)(\cstwo) \subset P_{f^t(c)}(0)$, and 
$J \cap \bigl(P_{f^t(c)}(0)-(f^t \circ g^k)(\cstwo)\bigr)$ is compactly
contained in $P_{f^t(c)}(0)$. \newcrit{And $f^t \circ g^k$ is \emph{univalent}
(in fact $f^t \circ g^k \on{\closure{\cstwo}}$ is an embedding).}

So, given any level $s$ Yoccoz piece $P(s)$, we have the branched
covering map $f^s: P(s) \> P_{f^t(c)}(0)$ for some $0 \le t < q$.
If $k$ is sufficiently large (given $P(s)$), then 
\mpar{Why take the closure here?}
$$
\{f^i(0) \mid 0 < i \le s\} \cap 
\closure{(f^t \circ g^k)(\cstwo)} \cap P_{f^t(c)}(0) = \nullset.
$$
\crit{could explain why I put ``$\cap P$'' in there}
Then for each point $z \in \partial P(s) \cap J$, $f^s(z) = \alpha$,
and we can define a single-valued branch of $f^{-s}$ on 
$(f^t \circ g^k)(\cstwo))$ such that $f^{-s}(\alpha)=z$, and
$f^{-s}((f^t \circ g^k)(\cstwo)) \subset P(s)$. 
Then we define the embedding $f_z: \closure S \> \closure{ P(s)}$
by $f_z=f^{-s} \circ f^t \circ g^k \circ \emb$. 
It is then readily seen
that $f_z(\closure S) \cap f_{z'}(\closure S) = \nullset$ for 
$z, z' \in \partial P(s)$, $z \neq z'$, and that  
$$
J\cap (P(s) - \bigcup_{z \in {\partial P(s)}} f_z(\closure S))
$$
is compactly contained in $P(s)$. Of course,
$f_z \on {S -\closure N}$ is quasiconformal, 
and $f_z(\closure S - \closure N) \subset \cx -J$, because
$\emb$ has these properties, and $f_z$ is just $\emb$ followed by 
(positive and negative) powers of $f$. 
Thus we have verified Lemma \ref{getS} for an arbitrary 
level $s$ piece
$P(s)$.

\chapter{The Tiling Lemma}
\label{decomp chapter}
\def\nondeg{N}
\def\romulus{N_1}
\def\remus{N_2}
\def\founders{N_i}
\def\targetlevel{L}
\def\annuli{{\bf A}}

\def\desc{descendant}			
\def\critdesc{critical descendant}	
\def\critdescendant {\critdesc}
\def\nncdesc{not-necessarily-critical descendant}
\def\critchild{child}
\def\pet{univalent descendant}
\def\cousins{cousins}
\def\grandchild{cousin}
\def\grandchildren{cousins}
\def\critchildren{children}

All the numerical variables in this chapter (as opposed to
object variables, like puzzle pieces) will denote integers.

Recall from Chapter \ref{overview} the statement of  the Tiling Lemma,
 \ref{tiling}:

\smallskip
\noindent {\em \statetiling}
\smallskip

In this chapter we break down the proof of Lemma \ref{tiling} into
three mutually exclusive cases for $f$, and settle each one by choosing
the $Q_i$ by a ``greedy algorithm'', and setting $R$ to the leftover
portion of the Julia set.  The first case, in which $\exists n:f^n(0)
=\alpha$, is trivial: we let $\{Q_i\}=\{P\}$, and $R=\nullset$. This case
must be eliminated in order to properly discuss the other two cases.
In the second case, the critically non-recurrent case, the leftover
set $R$ comprises at most one point.  In the third case, the
critically recurrent case, the leftover set $R$ is a Cantor set.

\section{List of cases}
Here are the three cases:
\begin{enumerate}
\item Some iterate of the critical point lands on the internal fixed point: 
$\exists n:f^n(0) =\alpha$.
\item The critical point is non-recurrent
($0 \notin \closure{\{f^n(0) \mid n > 0\}}$), but case 1 does not hold
($\forall n, f^n(0) \neq \alpha$).
\item The critical point is recurrent:  $0 \in \closure{\{f^n(0) \mid n > 0\}}$.
\end{enumerate}
We quickly treat case 1 in section \ref{trivial}.  After introducing
some notation in section \ref{notes} for cases 2 and 3, and making
some basic observations, we take care of case 2 in section
\ref{cnr}. We set up the proof for case 3 in section \ref{cr} and
finish it in section \ref{ver}.
\section{Proof for the $\exists n:f^n(0) =\alpha$ case} \label{trivial}
 In this case, $0 \in \Gamma_m$ for all $m \ge n$, so we can let
$L=n$, and given any piece $P$ of level $m>n$, $f^{m-n}\on P$ is
univalent, so we can form the trivial decomposition of $P$, namely
$T=R=\nullset$, and $Q_1=P$.

\section{Notation and setup, assuming $\forall n, f^n(0) \neq \alpha$}
\label{notes}

For $z \in J$, denote by $\piece n z$ the level $n$ puzzle piece that
contains $z$.  This is well-defined if $f^n(z) \neq \alpha$, but if
$f^n(z) = \alpha$, there will be no level $n$ piece containing $z$,
and more than one level $n$ piece containing $z$ in its closure.  So
this notation can only be used when we know already that $f^n(z) \neq
\alpha$.  In particular, our assumption here that $\forall n, f^n(0)
\neq \alpha$, ensures that $\piece n 0$ will be defined for all $n$.
We will call $\piece n 0 $ the {\em critical piece\/} of level $n$.

Given $z,n$ such that $z\in J$ and $f^n(z) \neq \alpha$, let $\annulus
n z$ denote $\piece n z - \closure{\piece {n+1} z}$.
 We call such a $\annulus n z$ a \term{combinatorial annulus} 
(even though it is not necessarily an annulus). 
\crit{Move definition of geometric annulus here?}  We will
call $\annulus n 0$ the {\em critical combinatorial annulus\/} of level $n$.

If we can prove \lt\ for the {\em critical\/} pieces of level greater
than $L$, then we can prove it for all pieces: given a piece $P$ of
level greater than $L$, it maps univalently by some iterate of $f$
either to some piece of level $L$ or to a critical piece of level
greater than $L$. In the former case, we can decompose $P$ trivially,
while in the latter case, we can pull back the decomposition of the
critical piece to $P$.

\section{Proof for the critically non-recurrent case, 
with $\forall n,$ $f^n(0) \neq \alpha$}
\label{cnr}


In this case the critical point forward orbit, $\{f^n(0)| n \in
\zplus\}$, is disjoint from some critical piece $\critpiece N$,
because the diameters of the $\critpiece N$ go to zero. We then set
$L = N$. Given a critical piece $P=\critpiece m$ with $m>N$, we let the
$\{Q_i\}$ be $\{\piece kz \mid k>m \text{ and }z \in \critannulus {k-1}\}$.  Thus
the the $Q_i$ are the level $k$ pieces that are subsets of $\critannulus {k-1}$.
Note that such a $Q_i=\piece kz$ is a subset of $\critpiece {k-1}$, but is
not $\critpiece k$.  If $f^t(\piece kz)$ were critical, for $t<k-L$,
then $f^t(\critpiece {k-1})$ would be critical, in which case $f^t(0)
\in \critpiece {k-t-1}
\subset \critpiece L$, a contradiction. We let $R=\{0\}$, which is the
intersection of all the critical pieces (by Theorem \ref{dzero}).  Finally, let
$T=\critpiece m \setminus (R \union \bigcup
\closure{Q_i})$.
 The only property of the $T,R,Q_i$ left to verify is that $T \cap J
=\nullset$, which is equivalent to $P \cap J \subset R \union \bigcup
\closure{Q_i}$.

So note that $J \cap \bigcup_{z \in \critannulus {k-1}} \closure{\piece kz} =
	J\cap \closure{\critannulus {k-1}}$, 
so $J \cap \bigcup \closure{Q_i} 
= J \cap (\closure{\critpiece m} - \bigcap_{l > m} \critpiece l)$. 
 But by Theorem \ref{dzero} (diameters of pieces go to zero), 
$\bigcap_{l > m} \critpiece l=\{0\}$. 
Therefore $J \cap \bigcup \closure{Q_i} = J \cap (\closure{\critpiece m} - R)$, 
which is equivalent 
to  $P \cap J \subset R \union \bigcup \closure{Q_i}$.

 We have thus shown Lemma
\ref{tiling} in the case where $f$ is critically non-recurrent.

\section{Notation and setup for the critically recurrent case}
\label{cr}
 For the critically recurrent case of \lt, we will need to let $R$ be
substantially more complicated. We say that $E \subset \cx$ is {\em
well-surrounded\/} if $E$ is compact and there exists a collection
\a\ of disjoint annuli in $\cx-E$ such that,
 if $x \in E$, the sum of the moduli of the annuli
in \a\ that surround $x$ diverges.

\bp \label{wshr}
If $R \subset \cx$ is well-surrounded, then $R$ is holomorphically removable. 
\ep

\begin{proof}
\ignore{
First, we observe that the well-surroundness of $R$ is a conformal invariant of
$\cx-R$, that is, if $R'$ is a compact subset of $\cx$, and $h: \cx-R \> \cx-R'$
is a conformal isomorphism (that maps $\infty$ to $\infty$), then $R'$ is
well-surrounded if $R$ is well-surrounded. 
}

We say that  compact subset $S$ of $\cx$
has {\em absolute area zero \/} if, whenever $S'$ is a compact subset
of $\cx$, and $h: \cx-S \> \cx-S'$ is a conformal isomorphism (that
maps $\infty$ to $\infty$), then $S'$ has measure 0.
Then the proposition follows
immediately from the following two results:
\begin{enumerate}
\item
\begin{thm}[McMullen]
A well-surrounded set has absolute area zero.
\end{thm}
This appears as Theorem 2.16 in \cite{McMullen:book:CDR}.
\ignore{
See  \cite{Branner:Hubbard:cubicsII}; see also Lyubich \cite{Lyubich:area:zero}.
 A precursor to this result (indeed also to Proposition \ref{wshr}) appears in
\cite{Sario:Nakai:book}.
}
\item
\begin{thm}\label{aaz}
A set that has absolute area zero is holomorphically removable.
\end{thm}
See \cite{Ahlfors:Beurling:nullsets}.
\end{enumerate}
\end{proof}

We also need some more facts about the Yoccoz partition.
 Here, as always, we assume that $f$ is
not renormalizable. 

If $\closure{\piece {n+1} z} \subset \piece n z$, 
then $\annulus n z$ is a (geometric) annulus.
\crit{move earlier?}


The following two statements can be found in the expositions of Milnor
	\cite{Milnor:local:connectivity} and Hubbard
	\cite{Hubbard:local:connectivity}:

\bl  \label{nondegann}
 There exists an $n$ such that $\critannulus n$ is an annulus. \el

In this case we call $\critannulus n$ a \term{critical annulus}. 
We say that $\critannulus m$ is a \term{\critdescendant} of $\critannulus n$
 if $f^{m-n}$ maps $\critannulus m$ onto $\critannulus n$ as an
unramified cover. 
 Note that in this case, if $\critannulus n$ is a geometric annulus, then so is
$\critannulus m$. If $f^{m-n}$ has degree 2, we say that
$\critannulus m$ is a \term{child} of $\critannulus n$. 

\bp  \label{crit diverges}
If $f$ is critically recurrent, then the sum of the moduli of the
\critdesc s of any critical annulus $A_n$ diverges.
\ep

We will assume for the rest of this section\newcrit{chapter} that $f$ is
critically recurrent. Let us now prove
\lt\ for this case. \newcrit{Why say this now?} 

We call two critical descendants of the same critical annulus \term{fraternal}
if neither one is a descendant of the other. Note that in this case they have no
descendants in common \newcrit{explain?}. We will need the following
lemma:
\bl \label{grandkids}
Every critical annulus $\critannulus{n}$ has at least two fraternal descendants 
	$\critannulus{n_1}$ and $\critannulus{n_2}$.  
\el

\begin{proof}
First note that, if $\critannulus m$ is a descendant of $\critannulus n$, and
$\critannulus l$ is a descendant of $\critannulus m$, then 
$\critannulus l$ is a descendant of $\critannulus n$. It follows that
either the above lemma is true for a given critical annulus $\critannulus n$,
or the descendants of
$\critannulus n$ form a sequence
$\critannulus{k_1}, \critannulus{k_2},
\critannulus{k_3}, \dots$ such that $\critannulus{k_j}$ is a descendant of
$\critannulus{k_i}$ whenever $j>i$. (So, in particular, $\critannulus{k_{i+1}}$
is a descendant of $\critannulus{k_i}$). But since the modulus of a descendant of
$\critannulus{m}$ is at most half the modulus of $\critannulus m$, the sum of
the moduli of the descendants of $\critannulus n$ would in this case converge, a
contradiction of Lemma \ref{crit diverges}. 
\end{proof}

\ignore{

We say that $\critannulus m$ is a \term{\critchild}\mpar{put definition of
child with descendant?} of
$\critannulus n$ if 
 $m > n$ and $f^{m-n} : \critannulus m \rightarrow \critannulus n$ is a 
	degree two  covering map. Note that this occurs if and only if
$f^{m-n} : \critpiece m \rightarrow \critpiece n$ and 
$f^{m-n} : \critpiece {m+1} \rightarrow \critpiece {n+1}$ are degree two branched covers.

We will call a \critchild\ of a \critchild\ a \grandchild. \newcrit{So a
grandchild is a kind of descendant}.

\begin{proof}
We quote two results from Milnor \cite{Milnor:local:connectivity}:
\bl \label{keepline}
Every critical annulus has at least one \critchild.
\el
and
\bl  \label{onlychild}
If a critical annulus has only one \critchild, then that \critchild\ in turn 
has at least two \critchildren.
\el
\newcrit{Could give proofs of these, following Lyubich.}
(In fact, a stronger statement than \ref{onlychild} is proven in
\cite{Milnor:local:connectivity}: every descendant of that only \critchild\ has
at least two \critchildren.  From this stronger statement and Lemma
\ref{keepline} above, Proposition \ref{crit diverges} follows
readily.\newcrit{So why not give the proof?}) 
So either our given annulus has at least two children, and each of them has
at least one child, or it has only child, which in turn has at least two
children. In either case our given annulus has two grandchildren.
\end{proof}

}

So, let $\nondeg$ be the level of the non-degenerate critical annulus
given by Lemma \ref{nondegann}, and let $\critannulus \romulus$ and
$\critannulus \remus$ be two fraternal descendants of $\critannulus
\nondeg$; their existence is guaranteed by Lemma
$\ref{grandkids}$.
 Set\mpar{display for emphasis?}
$L=\max (\romulus, \remus)+3$.

Now let $P=\critpiece \piecelevel$, with $\piecelevel > L$. We choose
the $Q_i$ by a kind of ``greedy algorithm''. First consider the set of all
pieces
$Q$ contained in $P$ such that $f^{q-L} \on Q$ is univalent, where $q$
is the level of $Q$. Then let the $Q_i$ be those elements of this set
that are not a sub-piece of any other element. Then, by the Markov
property of the Yoccoz partition, the $Q_i$ are automatically mutually
disjoint, and between them they cover as much of $P$ as we could hope
to cover. We let $R=(P \setminus \bigcup \closure{Q_i}) \cap J$, and
let $T= P \setminus (R \cup \bigcup \closure{Q_i})$.
Now we just need to verify is that $R$ is compact and
holomorphically removable. We will do so by showing that $R$ is
well-surrounded. In order to do this, we must of course define a set
$\annuli$ of annuli.

We let 
$$
\annuli=\{ \annulus nz \mid z \in R, n \ge\piecelevel, \text{ and } 
f^{n-\nondeg}:
\annulus nz \> \critannulus {\nondeg}\text { is a covering map.}\}
$$
  (We will verify that $z \in R$ implies that $\forall n, f^n(z)
\neq \alpha$, so $\annulus nz$ is well-defined).

We now need only  the following:
\begin{lem} The set $R$ defined above is well-surrounded by $\annuli$. In
particular,
\begin{enumerate} 

\item \label{num compact}the set $R$ is compact, 
\item \label{num disjoint} the annuli in $\annuli$ are mutually disjoint, and
disjoint from $R$, and
\item \label{num diverges} the sum of the moduli of the annuli in $\annuli$
that surround any given point in $R$ diverges. 
\end{enumerate}
\end{lem}

To verify \ref{num compact} and \ref{num disjoint} we need just $L >
\nondeg$; it is for one case of the verification of \ref{num diverges} that 
we need $L >
\max{\founders +3}$.
\newcrit{What happened to the explanations?}

\section{Proof of well-surroundedness of $R$ for the critically recurrent
case}
\label{ver}


Here are the verifications of the above three statements.

\subsection{Compactness of $R$}
\def\P{P}
\begin{lem}\label{turd1}
If $\P$ is a piece, and $\eta \in \partial \P \cap J$, 
then $P$ has a subpiece $\P'$
of level $k$ with $\eta  \in \partial \P'$ and $f^{k-L}$ univalent on $\P'$.
\end{lem}

\begin{proof} Note that $\alpha \notin \closure {\critpiece {\nondeg + 1}}$, because
$\closure {\critpiece {\nondeg + 1}} \subset \critpiece \nondeg$, and
$\alpha \notin \critpiece \nondeg$. Therefore, no piece of level greater than
$\nondeg$ with $\alpha$ on its boundary can be critical. Now, for any piece $\piece
nz$, if $\alpha \in \partial \piece nz$, then 
$\alpha \in 
\partial (f^k(\piece nz))$, for all $0 <k \le n$.
 Therefore, if $z,n$ are such that $n > \nondeg$ and
 $\alpha \in \partial \piece nz$, then $f^{n-\nondeg}:\piece nz \> \piece
 \nondeg {f^{n-\nondeg}(z)}$ is univalent. In other words, every piece
 with $\alpha$ on its boundary maps univalently (by an iterate of $f$)
 to a level $\nondeg$ piece (and, hence, to a level $L$ piece, since $L > N$).

Now, given $\eta \in \partial \P \cap J$, chose $m$ such that $f^m(\eta) = \alpha$.
 Note that $\eta$ is not a critical point of $f^m$,
 due to our standing assumption in this section that 
$\forall n>0, f^n(0) \neq \alpha$. Therefore, there is some level $k_0 \ge m$ such that
$f^m$ is univalent on  every piece $\piece kz$ of level $k \ge k_0$  with $\eta \in
\partial \piece kz$.  But then $\alpha \in f^m(\piece kz)$, so if $k-m > L$, then 
$f^{(k-m)-L} \on {f^m(\piece kz)}$ is univalent, so, in any case, $f^{k-L}\on
{\piece kz}$ is univalent. So the desired subpiece $\P'$ is the unique piece of
level $k$ (say with $k = k_0$) that is contained in $P$ and has $\eta$ on its
boundary. 
\end{proof}

Note that the $P'$ described above will be contained in some $Q_i$, because it
maps univalently up to the level $L$. 
Note also that $P'$ contains the intersection of $P$ with some neighborhood
of $\eta$.

\ignore{
Thus $\eta \in \closure{Q_i}$, so $\eta
\notin R$. We conclude that no point of $R$ is on the boundary of any piece, so
that $\piece nz$ is well-defined for $z \in R$. 
}

\begin{cor}
If $x_i \rightarrow x$, and $x_i \in R$, then $x$ is not the boundary of any
piece.
\end{cor}
For if $x$ were on the boundary of some piece, then we could choose a
subsequence of the $x_i$ to lie in one of the finitely many pieces of
a given level that have $x$ as a boundary point. But then the
preceeding lemma provides a contradiction, because no points in that
piece that are sufficiently close to $x$ can be in $R$.

\begin{cor}
If $\eta \in \partial \P$ for some piece $\P$, then $\eta \notin R$. 
\end{cor}
This is because $\eta \in \closure{\P'} \subset \closure Q_i$ for some $Q_i$.
\begin{cor} 
If $z \in R$, then $z$ is not on the boundary of any piece, so $\piece
zn$ (and $\annulus zn$) is well-defined.
\end{cor}
This is an immediate consequence of the previous corollary.

\begin{lem}
If $x_i \in R$, and $x_i \rightarrow x$, and $x$ is not on the boundary of a piece,
then $x \in R$.
\end{lem}

\begin{proof} If $x$ were not in $R$, 
then there would be a piece  containing it that
maps univalently up to level $L$. But then a whole neighborhood of $x$
would not be in R.
\end{proof}

We conclude from the above that $R$, in any piece, is compact (and in
particular stays away from the boundary of that piece). We also conclude that
$\piece nz$ (and $\annulus nz$) is well-defined for all $n \ge 0$ and $z \in R$.

\subsection{Disjointness of Annuli}

\begin{lem} \label{burp}
No annulus in $\annuli$ can contain a point of $R$ in its closure.
\end{lem}
{
\def\A{\annulus kz}
\begin{proof} 
Every annulus $\annulus kz$ in $\annuli$ is composed of pieces of level greater
than $\piecelevel$ that map univalently to level $\nondeg$ (because $\A$ is an
unramified cover of $\critannulus \nondeg$, and an iterate of $f$, restricted to
any piece, either has a critical point or is univalent), and we assume that $L
\ge \nondeg$+1. So every point in $\closure \A$ lies in some $\closure{Q_i}$, and
hence cannot be in $R$.
\end{proof}
}

\begin{lem}\label{intersect}
Suppose two combinatorial annuli, $\annulus kz, \annulus lw$, intersect. Then 
\begin{enumerate}
\item $z \in  \closure{\annulus lw}$, or
\item $w \in \closure{\annulus kz}$, or
\item $\annulus kz = \annulus lw$.
\end{enumerate} 
\end{lem}

\def\W_#1{\piece {#1}z}
\def\V_#1{\piece {#1}w}
\begin{proof}
Recall that $\annulus kz = \piece kz -\closure{\piece {k+1}z}$, and
 $\annulus lw = \piece lw -\closure{\piece {l+1}w}$. If $k=l$, then $\piece kz =
 \piece lw$, and either $\piece {k+1}z = \piece {l+1}w$ (case 3), or $\piece
 {k+1}z \cap \piece {l+1}z = \nullset$. If the latter holds,  then we have (since
 pieces are open)
 $\piece {l+1}w \subset \piece kz - \closure {\piece{k+1}z}$, and then
 $w \in \piece{l+1}w \subset \annulus kz$.  If $k > l$, then $\piece
 kz \subset \piece lw$ (because $\piece
 kz \cap \piece lw \neq \nullset$),
 and $\piece kz \not\subset \piece {l+1}w$, so
 in fact $\W_k \cap \closure{\V_{l+1}} = \nullset$  by the
Markov property for pieces, so $z \in \W_k \subset
\V_l - \closure{\V_{l+1}}=\annulus lw$.
\end{proof}

\begin{cor}
No two distinct annuli in $\annuli$ can intersect.
\end{cor}

\begin{proof}
 Suppose there were two distinct annuli $\annulus kz, \annulus lw \in \annuli$,
 with $z, w \in R$, and $\annulus kz \cap \annulus lw \neq \nullset$.
Then by \ref{intersect}, either $z \in  \closure{\annulus lw}$, or
$w \in \closure{\annulus kz}$. But this contradicts Lemma \ref{burp}.
\end{proof}

\subsection{Divergence}

\def\seq#1#2#3{(#1)_{#2}^{#3}}
\def\sseq#1{(#1)}
\def\passes{rises past}
\def\taulike{rise-and-drop}
\def\a#1{a_{#1}}
\def\step{step}
\def\modulus{\operatorname{\text{mod}}}

For\mpar{We want to {\em improve} this paragraph---be more explicit about its
direction, if possible.}
 $z \in J$, let $\annuli_z$ denote all elements of the form
$\annulus nz$ of $\annuli$, that is, all elements of $\annuli$ that
surround $z$.  Then our goal is to show, for each $z \in R$, that the
sum of the moduli of the elements of $\annuli_z$ diverges. The first
step is to determine which $n$ are such that $\annulus nz \in
\annuli_z$.
 This is done with the aid of the function $\taufunc nz$, first
defined by Shishikura (following Yoccoz) in his proof of Theorem \ref{mzero}.  
Then one property of
$\taufunc nz$ is abstracted in \taulike\ functions, defined below.  We
prove certain lemmas about
\taulike\ functions.  One such lemma is enough to deduce divergence of
the sum of the moduli of the elements of $\annuli_z$ in the case where
$\sup \taufunc nz$ is infinite. In this case divergence is deduced from
the divergence of $\annuli_0$, quoted as Lemma \ref{crit diverges} in
the previous section. In the other case (when $\sup \taufunc nz$ is finite), we
show that $\annuli_z$ contains infinitely many copies of one of the two
$\critannulus{\founders}$, and hence the the sum of the moduli of the annuli in
$\annuli_z$ diverges.

Given $n \in \nats,z \in J$ such that $f^n(z) \neq 0$, there is at most one $m \in [0,n]$ such that
$f^{n-m}(\piece n
	z) = \critpiece m$, and $f^{n-m} \on {\piece nz}$ is univalent
	(so then $f^{n-m} : \piece nz \> \critpiece m$ is an
	isomorphism).\newcrit{Should I explain $f^0$?}
If such an $m$ exists, then we set $\taufunc nz = m$.  If no such $m$
exists, then $f^n\on {\piece nz}$ is univalent, and $f^n (\piece nz)$
is not a critical piece.  In this case we set $\taufunc nz = -1$.

So now we can write
$$
R = \{z \in J \cap \critpiece \piecelevel \mid 
\forall n \ge \piecelevel, \quad
f^n(z) \neq \alpha \text{ and }\taufunc nz >L  \}.
$$

We will be interested in the values of $\taufunc nz$ for $z \in R$ and
$n \ge \piecelevel$. In particular, by our definition of $R$, $\taufunc 
nz$ will be non-negative. In the statements that follow, we will have
the standing assumption that $\taufunc nz$ is non-negative, whenever
$n,z$ are mentioned in the hypothesis.

\begin{lem} \label{get}
If $\taufunc n z = m \ge 0$, and $\taufunc {n+1} z = m+1$, then $f^{n-m}:
\annulus nz \to \critannulus m$ is an isomorphism.
\end{lem}
\begin{proof}
 We have that $f^{m-n} :
 \piece n z \to \critpiece m$ is an isomorphism,
 and $f^{m-n}(\piece {n+1} z) = \critpiece {m+1}$, so
 $f^{m-n}(\piece n z - \closure{\piece {n+1} z}) =
\critpiece m - \closure{\critpiece {m+1}}$. 
\end{proof}
\begin{cor} \label{equal modulus}
If $\critannulus m$ is a descendant of $\critannulus \nondeg$, and 
$\taufunc n z = m \ge 0$, and $\taufunc {n+1} z = m+1$ (for $n \ge p$), then 
$\annulus nz \in \annuli_z$  (and the modulus of $\annulus nz$ is equal to the 
modulus of $\critannulus m$).
\end{cor}
\begin{proof}
By Lemma \ref{get}, $f^{n-m}: \annulus nz \to \critannulus m$ is an isomorphism.
By assumption, $f^{m-\nondeg}:\critannulus m \> \critannulus \nondeg$ is a
covering map. Therefore $f^{n-\nondeg}:\annulus nz \> \critannulus \nondeg$ is a
covering map, so $\annulus nz \in \annuli_z$. 
\end{proof}

We now make a simple observation about the function $\taufunc nz$:
\begin{lem}
\label{tau taulike} For all $n,z$, $\taufunc {n+1} z \le \taufunc n z +1$
\end{lem}
\begin{proof}
For all $\nu,\zeta$, we have that $f^{\nu- \taufunc \nu \zeta}$ is the
greatest iterate of $f$ that is univalent on $\piece \nu \zeta$.
Therefore $f^{n - \taufunc n z} \on {\piece n z}$ is univalent, and
$\piece {n+1} z \subset \piece n z$, so $f^{(n+1) - (\taufunc n z+1)}
\on {\piece {n+1} z}$ is univalent, and therefore $n+1 - \taufunc {n+1} z \ge n -
\taufunc n z$, that is, $\taufunc {n+1} z \le \taufunc n z +1$.
\end{proof}

This then motivates the following definition:
\begin{defn}
A sequence of non-negative integers $\sseq {\a n}$ is
\term{\taulike} if it is bounded
below, and
$\forall n,
\a {n+1}
\le
\a n +1$. The sequences we will consider will either be finite in length or
forward-infinite.
\end{defn}
So, by Lemma \ref{tau taulike}, $\taufunc nz$ is rise-and-drop for all $z \in
J$. \crit{What about when $\tau nz=-1$?}

\begin{defn}
A \term{\step} is a pair $(m, m+1)$ of consecutive non-negative integers.
\end{defn}
\begin{defn}
We say that a \taulike\ sequence $\sseq {\a n}$ \term{\passes} a step
$(m, m+1)$ at time $(n, n+1)$ if $\a n = m$ and $\a {n+1} = m+1$
\end{defn}
Note that if $\critannulus m$ is a descendant of $\critannulus \nondeg$,
and $\taufunc nz$ rises past
$(m, m+1)$, at time
$(n, n+1)$, then
$\annulus nz \in \annuli_z$, and
$\modulus \annulus nz = \modulus \critannulus m$.  

\bl[Intermediate value theorem for \taulike\ sequences] \label{IVT}\quad
Suppose that $\seq {\a i}{i=1}l$ is \taulike. Then if 
$k \le l$ and $a_k \le m < m+1 \le a_l$,  then $\seq {\a i}{i=1}l$ rises
past $(m, m+1)$. 
\el
\ignore{
\begin{proof}
Let $t= \inf \{i \mid \a i \ge m+1\}$. 
Then $t > k$, $\a t \ge m+1$, and $\a {t-1}  < m+1$, so $\a {t-1} \le m$, and
since $\a t \le \a {t-1} +1$, we must have $(\a {t-1}, \a t) = (m, m+1)$. 
\end{proof}
}
\begin{proof}
Let $s=\sup \{i \mid \a i \le m\}$. Then $\a s \le m$, $\a {s+1} \ge
m+1$, and $\a {s+1} \le \a s +1$, so $(\a s, \a {s+1}) = (m, m+1)$.
\end{proof}

We now present two further lemmas on \taulike\ sequences. The first is for the
case where $\sup \taufunc nz = \infty$, and the second is for the case where
$\sup \taufunc nz$ is finite. 

\bl \label{unbounded}
If
$\seq{\a n}{n=k}{\infty}$ is \taulike, and $\sup \seq{\a
n}{n=k}{\infty}$ is infinite, then $\sseq{\a n}$
\passes\ all but finitely many \step s.
\el
\def\infinum{b}
\begin{proof}
The given sequence $\a n$ is bounded below, so let $\infinum = \inf \a
n$.  We will show that $\a n$ rises past $(m, m+1)$ for all $m \ge
\infinum$.  Since $\a n$ is discrete-valued, its infinum is realized,
so let $k$ be such that $\a k = \infinum$.  Now, suppose we are given
such an $m$.  Then, since the sequence $\sseq {\a i}$ has a discrete domain,
$\lim \sup \a i = \sup \a i$, so $\exists l > k$ such that $\a l \ge
m+1$.  Then, by our ``intermediate value theorem''(\ref{IVT}),
$\exists s$, with  $k \le s < s+1 \le l$, such that $(\a s, \a {s+1}) = (m,
m+1)$.
\end{proof}

\bl \label{bounded}
Suppose
$\seq{\a n}{n=k}{\infty}$ is \taulike, and $\sup \seq{\a
n}{n=k}{\infty}$ is finite.  Then 
\begin{enumerate}
\item 
$\sseq{\a n}$ makes the same drop
infinitely often: $\exists r, s$ with $r \ge s$ such that $\a n = r$
and $\a {n+1} =s$ for infinitely many n. 
\item 
If  $m$ is  given such that $s \le m < m+1 \le r$,
then $\sseq {\a n}$ \passes\ $(m, m+1)$ infinitely many times.
\end{enumerate}
\el
\begin{proof}
In this case, $\a n$ realizes only finitely many values, so there are
 only finitely many possible pairs of values $(\a n, \a {n+1})$ with
 $\a {n+1} \le \a n$, so at least one such pair of values must be
 realized infinitely often.  So there is a monotonically increasing
 sequence $\seq {n_i}{i=1}{\infty}$, with $\a {n_i} =r$ and $\a {n_i
 +1} =s$ for some $r \le s$.  Then, given $m$ with $r \le m < m+1 \le
 s$, we note that, for each $i \in \nats$, $\a {n_i+1} =s$ and $\a
 {n_{i+1}} =r$, so, by our ``intermediate value theorem''(\ref{IVT}),
 there exist $s_i$, with $n_i +1 \le s_i < s_i+1 \le n_{i+1}$, such
 that $(\a {s_i}, \a {s_i +1}) = (m, m+1)$.
\end{proof}

With the help of Lemma \ref{unbounded}, we can now settle the case where $\sup
\taufunc nz = \infty$.
\begin{lem} \label{one}
If $\sup \taufunc n z= \infty$, and $z \in R$, then the sum of the
moduli of the annuli in $\annuli_z$ diverges.
\end{lem}
\begin{proof}
By Lemma \ref{tau taulike}, $\seq {\taufunc n z}
 {n=\piecelevel}{\infty}$ is \taulike, so by Lemma \ref{unbounded}, it
 rises past all but finitely many \step s. By 
Corollary \ref{equal modulus}, for
 each time it \passes\ a \step\ $(m, m+1)$ with $\critannulus m$ a
 descendant of $\critannulus \nondeg$, we get an element of $\annuli_z$
 with modulus equal to the modulus of $\critannulus m$. Therefore, by
 Lemma \ref{crit diverges}, the sum of the moduli of the annuli in
 $\annuli_z$ diverges.
\end{proof}

\medskip

The rest of this subsection is devoted to showing that the sum of the
moduli of the annuli in $\annuli_z$ for $z \in R$ diverges
when $\sup \taufunc n z$ is finite. In this case the forward orbit of
$z$ does not accumulate on the critical point, and we cannot pull back
a copy of each of the annuli surrounding the critical point.  Instead,
we will show that every time $\taufunc n z$ fails to increase by 1
(when $n$ increases by 1), it in fact ``drops'' past a
\step\ corresponding to the level of a \critdesc\ of $\nondeg$.
 The argument here is a little technical: we in fact show that it drops
past the level of one of the two fraternal descendants $\critannulus \founders$ of
$\critannulus \nondeg$.

Then we can conclude, using Lemma \ref{bounded} and Lemma \ref{get} (or Corollary
\ref{equal modulus}), that $\annuli_z$ contains infinitely many conformal copies of
a single
\critdesc\ of
$\critannulus
\nondeg$. So in this case the series of moduli of
elements of $\annuli_z$ contains infinitely many
copies of the same number, and therefore diverges.

\newcrit{What's the step-by-step game plan?}

\begin{lem} \label{firstborn}
Suppose $\critannulus n$ is critical annulus. Because the critical
point $0$ of $f$ is recurrent, there is some $m > 0$ such that $f^m(0) \in
\critpiece {n+1}$; chose the least such $m$.  Then $f^m: \critannulus {n+m} \to
\critannulus n$ is a double cover, so $\critannulus {n+m}$ is a child
of $\critannulus n$.
\end{lem}
\begin{proof}
We have $f^m(\critpiece {m+n+1})=\critpiece {n+1}$, so $f^m(\critpiece
{m+n})=\critpiece n$, and $f^m \on {\critpiece {m+n+1}}$ is a degree 2
branched cover, so all we need check is that $f^m \on {\critpiece
{m+n}}$ is degree 2.  If not, then $f^i(\critpiece {m+n})=\critpiece
{m+n-i}$ for some $0 < i < m$. But then $m+n-i \ge n+1$, so we get
$f^i(0) \in \critpiece {n+1}$, a contradiction.
\end{proof}
\begin{cor}\label{next}
If $\critannulus n$ is a critical annulus, and $f^k(0) \in \critpiece
 {n+1}$, then $\critannulus n$ has a child $\critannulus t$, with $t
 \le n+k$.
\end{cor}
\begin{proof}
 Use Lemma \ref{firstborn}: the $m$ in Lemma \ref{firstborn} satisfies
$m \le k$, so the child $\critannulus {m+n}$ satisfies $n+m \le n+k$.
\end{proof}

\begin{lem}  \label{getdesc}
Suppose $f^k(\critpiece {n+k}) = \critpiece n$. Then for all $l<n$
\newcrit{$l<n+k$ is enough},
 there exists $t$ such that $n+k > t \ge n$ and $\critannulus t$ is
 a descendant of $\critannulus l$.
\end{lem}

\begin{proof}
 Let $m$ be the greatest integer such that $m < n$, and $\critannulus
m$ is a descendant of $\critannulus l$. (We allow the possibility that
$m=l$.)  Then $f^k(0) \in \critpiece n \subseteq \critpiece {m+1}$, so
by Corollary $\ref{next}$, $\critannulus m$ has a child $\critannulus
t$ with $t
\le m+k < n+k$, but $t \ge n$ by our choice of $m$.
 So, in summary, $\critannulus t$ is a descendant of $\critannulus l$
(via $\critannulus m$), and $n+k > t \ge n$.
\end{proof}

\begin{lem} \label{drop nature}
If $\taufunc {n+1} z \le \taufunc n z$,  then 
$$
f^{\taufunc nz - (\taufunc {n+1}z -1)}(\critpiece {\taufunc nz})=
\critpiece {\taufunc {n+1}z -1}.
$$
\end{lem}

\begin{proof}
Let $a:=\taufunc {n+1} z$ and $b:=\taufunc n z$. 
Then $f^{n+1-a}(\piece {n+1}z)=\critpiece a$, so
$f^{n+1-a}(\piece {n}z)=\critpiece {a-1}$.
Also $f^{n-b} (\piece nz) = \critpiece b$. 
Therefore $f^{(n+1-a)-(n-b)}(\critpiece b)=\critpiece {a-1}$.
\end{proof}

\begin{lem} \label{findann}
Suppose $z \in R$, $a:=\taufunc {n+1} z \leq b:=\taufunc n z$, and $a > L$. 
Then there exists $m$ such that $a \leq m < m+1 \le b $,
and $\critannulus m$ is a \critdesc\ of one of the two $\critannulus \founders$. 
\end{lem}

\begin{proof}
We have $f^{b-(a-1)}(\critpiece b)=\critpiece a-1$ by Lemma  \ref{drop nature}. 
Then by Lemma \ref{getdesc}, since both of the
levels $\founders$ are less than $a-1$ (by the definitions of $L$ and $R$),
there is a descendant $\critannulus t$ of each $\critannulus \founders$ with $a-1 \le t
< t+1 \le b$. The two descendants are distinct, so one of the $t$'s
must satisfy $a \le t < t+1 \le b$. This $t$ is the required $m$.
\end{proof}

\begin{lem} \label{two}
If $\limsup \taufunc n z < \infty$, and $z \in R$, then the sum of the
moduli of the annuli in $\annuli$ that surround $z$ diverges.
\end{lem}

\begin{proof} 
\mpar{Need new letter for `$p$' here}
By the first part of Lemma \ref{bounded}, there exists $a \le b$
such that for infinitely many $n$, $\taufunc {n+1} z = a$ and
$\taufunc n z = b$.  Then by Lemma \ref{findann} we can find a
descendant $\critannulus m$ with $a \le m < m+1 \le b$.  So then by
the second part of Lemma \ref{bounded} there are infinitely many $q$
with $\taufunc q z = m$ and $\taufunc {q+1} z = m+1$.  Then by Corollary 
\ref{equal modulus}, $\annulus q z \in \annuli$, and $\modulus \annulus qz =
\modulus \critannulus m$. Thus there are infinitely many
annuli in $\annuli_z$ with modulus $\modulus
\critannulus m$.   So the sum of the
moduli of the annuli in $\annuli_z$ diverges.
\end{proof}

\begin{lem} \label{diverges}
For all $z \in R$, the sum of the moduli of the annuli in $\annuli_z$ diverges.
\end{lem}

\begin{proof}
This is just the conjunction of Lemmas \ref{one} and \ref{two}.
\end{proof}

\chapter{Further Results}  
\label{further}

\section{Local Connectivity of Corresponding Points in the Mandelbrot Set}
\label{MLC}
\def\state#1{{\noindent \sl #1. }}

In this section we will prove:
\begin{thm}
If\mpar{What about the finitely renormalizable case?} 
$c \in \partial M$, and $f_c(z)=z^2+c$ is not renormalizable and
has no indifferent fixed point, then $M$ is locally connected at $c$. 
\end{thm}
\newcrit{Must also cite G. Levin's result (and, of course, Yoccoz's original result}

The proof is by analyzing the behaviour of the graphs $\Gamma_n(c)$ 
as c varies. (Recall the definition of $\Gamma_n$ (here written 
as $\Gamma_n(c)$ to emphasize its dependence on $f_c$) from 
section \ref{Yoccoz Partition}.)

\begin{proof}
We have $f_c=z^2+c$\newcrit{cut this?}.  For $c \in M$ let $\beta(c)$ denote the landing
point of the zero ray, and for $c \not= 1/4$ let $\alpha(c)$ denote
the other fixed point.  Then if $\alpha(c)$ is repelling let
$\Gamma_n(c)$ denote the level $n$ Yoccoz graph for $f_c$.  We wish to
show that $M$ is locally connected at $c$ if $c$ is
non-renormalizable.

\medskip
The set of all $c$ for which $\alpha(c)$ is repelling and has rotation
	number $p/q$ is called the $p/q$ limb of the Mandelbrot set,
	denoted $M_{p/q}$.  There is a unique $c_{p/q}$ for which
	$f_c$ has a parabolic fixed point of multiplier $e^{2 \pi i
	p/q}$, and $M_{p/q}$ is one of the two components of
	$M-\{c_{p/q}\}$.

\begin{thm}
The\mpar{Need a reference---cf. Hubbard}
 arguments of the rays landing at a repelling periodic point of a polynomial
 (with connected Julia set) are stable under
perturbation. Likewise, the arguments of rays landing at a
non-ramified preiterate of a repelling periodic point are stable under
perturbation.  A non-ramified preiterate of $\alpha$ is a point $z$
such that $f^n(z)=\alpha$, and $(f^n)'(z) \not= 0$.
\end{thm}

Suppose $c \in M_{p/q}$.  Then we can define $\Gamma_n(c)$ for all
$n$.  For a given $n > 0$, if $f_c^n(0) \not= \alpha$ (here 0 is, of course, the
critical point of $f_c$), then $\Gamma_n(c)$ remains constant on some neighborhood
of $c$ in $M$, in the sense that the information of which rays land in groups of
$q$ at points in
$f^{-n}(\alpha)$\newcrit{Which groups of rays land together (cf. $\simeq$ relation
in Chapter 2)} remains constant in that neighborhood.  In other words, for any given
$n$, the combinatorial information in
$\Gamma_n(c)$ is locally constant in $M_{p/q} - \{ c \mid f^n_c(0) =
\alpha(c)\}$.  (Note that $\{ c \mid f^n_c(0) = \alpha(c)\}$ is finite, since it is
a subset of $\{c \mid f_c^n(0) = f_c^{n+1}(0)
\}$).  Therefore the information is constant on the finitely many components of
$M_{p/q} - \{ c \mid f^n_c(0) = \alpha(c)\}$, which are each open in $M$.

Now, 
	suppose $c \in M_{p/q}$, 
	and $\forall n, \, f_c^n(0) \not= \alpha(c)$.
\def\comp#1{M^{#1}_{p/q}}
Then, for each $n > 0$, consider $\comp n(c)$, the component of
	$M_{p/q} - \{ c \mid f^n_c(0) = \alpha(c)\}$ that contains
	$c$.  We claim that, if $f_c$ is non-renormalizable, then the
	sets $\comp n(c)$ form a neighborhood base for $c$ in $M$.

\def\closure{\overline}
\newcrit{This is confused. $\bigcap_{n=1}^\infty \closure{\comp n(c)} \neq \{c\}
\rightarrow \exists c', c' \neq c$ s.t. $c' \in \comp n(c) \forall n$ (because
continua are uncountable).} Consider the continuum $\bigcap_{n=1}^\infty
\closure{\comp n(c)}$. We just need to show that it is degenerate (i.e. a single
point).  If it is non-degenerate, then we can find $c' \not= c$ such that $\forall
n, \, f_{c'}^n(0) \not= \alpha(c')$ (and such that $c' \not=
c_{p/q}$).  Then $c' \in \comp n(c)$ for all $n$, so $\Gamma_n(c) =
\Gamma_n(c')$ for all $n$.  But then $c=c'$ by the following:

\begin{thm}
Suppose $c, c' \in M_{p/q}$ for some $0 < p/q <1$ in lowest terms.  If
$f_c$ and $f_{c'}$ are combinatorially equivalent (i.e. $\forall n, \,
\Gamma_n(c)$ and $\Gamma_n(c')$ have the same rays landing in groups
of $q$) and non-renormalizable, then $f_c=f_{c'}$.
\end{thm}
\begin{proof}

\state{Step 1}
If $c$, $c'$ are combinatorially equivalent and non-renormalizable,
	then there exists a homeomorphism $h: \cx \rightarrow \cx$
	such that $h \circ f_c=f_{c'} \circ h$ on $\cx$,
	$h(J_c)=J_{c'}$, and $h|_{\cx-J_c}$ is conformal.

\state{Step 2}
From\mpar{Move to after proof of Step 1?}
 the above and the holomorphic removability of $J_c$ we can
	immediately conclude that $h$ is conformal, and hence $f_c$
	and $f_{c'}$ are conformally and thus affinely conjugate, so
	$c=c'$.

\smallskip
\state{Proof of Step 1}
For all $n$, since $\Gamma_n(c)$ and $\Gamma_n(c')$ have the same
	combinatorics, there is a canonical homeomorphism from the one
	to the other (off of $\alpha(c)$ and its preiterates, it
	factors through the two Riemann maps).  That homeomorphism can
	be extended conformally outside of $\Gamma_n(c)$, and
	arbitrarily on the bounded components of the complement, to
	form a homeomorphism $h_n$ from $\cx$ to $\cx$.  The $h_n$ are
	eventually constant on the complement of $J_c$ (and converge
	uniformly on any compact subset of the complement), and are
	uniformly bounded on $J_c$.  Any pointwise limit\newcrit{really should use
Arzela-Ascoli instead} 
$h_\infty$ of
	the $h_n$ is conformal off of $J_c$, and, because the diameter
	of the pieces of $\Gamma_n(c)$ and $\Gamma_n(c')$ go to zero
	as $n \rightarrow \infty$, $h_\infty$ is continuous
	(sufficiently nearby points in the domain lie in the same or
	adjacent small pieces, and therefore have nearby images),
	injective (distinct points eventually lie in distinct and
	non-adjacent pieces, and hence the lim inf for the distance
	between their images under the $h_n$ is positive), and proper,
	and therefore it is a homeomorphism.  It is a conjugacy off of
	$J_c$, and therefore on all of $\cx$.
\end{proof}

\end{proof}

\section{Finitely Renormalizable Quadratic Polynomials}
\label{fr}
\def\R{{\cal R}}

\subsection{Definitions}\label{renorm def}
Suppose $f_c(z)=z^2+c$ has both fixed points repelling. Then we can form the
Yoccoz graph $\Gamma_n$. for $f$. Suppose futher that $f^n(0) \neq \alpha$ for all
$n > 0$. Then $\critpiece n$ is well defined. 

We say that $f$ is \term{combinatorially renormalizable} (with period $n>1$) if
$\exists k,n$ such that $f^n:\critpiece {k+n} \> \critpiece k$ is a degree two
branched cover, and
$f^{tn}(0) \in \critpiece {k+n}$ for all $t > 0$. We call the map
$f^n:\critpiece {k+n} \> \critpiece k$ a \term{combinatorial renormalization} (with
period $n$).
If such an $n$ exists for $f$, it
will be unique, and in fact, for all $k' >k$, $f^n:\critpiece {k'+n} \>
\critpiece k'$ will have the same properties (mentioned above) as $f^n:\critpiece
{k+n} \> \critpiece k$. 
The set $K_{\R f} := \{z \in \critpiece {k+n} \mid \forall t>0, f^{nt}(z) \in
\closure{\critpiece {k+n}}\}
=\bigcup_i\closure{\critpiece i}$ is called the filled-in Julia set of the combinatorial renormalization for $f$.
$K_{\R f}$ does not depend on the choice of $k$, so it is a well-defined object
(given a renormalizable map $f$). \newcrit{proof?}

Following Douady and Hubbard, we define a \term{quadratic-like map} as a
holomorphic degree 2 branched cover $g: U' \> U$, where $U', U \subset \cx$ are
topological disks, and $\closure{U'} \subset U$.  We define the filled-in
Julia set $K_g$ of $g$ by $K_g = \{ z \in U' \mid \forall t>0, g^t(z) \in U'$. We
say that  $g$ is \term{non-trivial} if the critical point of $g$ lies in $K_g$.
In this case, $K_g$ is connected. Given $f_c(z) = z^2 +c$, we say that $f$ is
\term{geometrically renormalizable} (with period $n$) 
if there exist $U', U$ and $n$ such that 
$f^n:U' \> U$ is a non-trivial polynomial-like map, and $0 \in U'$. 
(Note then that 0 is the unique critical point of $f^n$ in $U'$.) 

We have the following theorem\cite{Milnor:local:connectivity,
Hubbard:local:connectivity}, which is part of the Yoccoz theory:
\begin{thm}[Straightening Theorem] \label{rcg}
If $f$ is combinatorially renormalizable with period $n$, then $f$ is
geometrically renormalizable with period $n$, and the Julia set of the
combinatorial renormalization is the same as that of the geometric renormalization.
\end{thm}
The converse is also true (but we will not need it), if we assume that $f$ is
\term{simply} (geometrically) renormalizable. For a definition of simple
renormalization, and a discussion, see
\cite{McMullen:book:CDR}. The geometric renormalizations that arise from the above
theorem will always be simple renormalizations. 

We will use the term \term{renormalizable} to mean combinatorially renormalizable,
which, by the above, is equivalent to being (simply) geometrically renormalizable. 

We will require the following theorem of Douady and Hubbard\cite{dh:plm}:
\begin{thm}
If $g: U' \> U$ is a quadratic-like map, then there exists a quasiconformal
embedding $h: U \> \cx$ and a map $f_c(z) = z^2 +c$ such that $h(g(w))=f_c(h(w))$
for all $w \in U'$. It follows then that $h(K_g) = K_f$. Moreover, $c$ is unique
if we require that $g$ is non-trivial, and that the dilatation of $h$ be zero
a.e. on $K_g$. 
\end{thm}
(Note that the last condition is trivially satisfied if $K_g$ has measure 0.)
In the case where $g$ is a (geometric) renormalization of $f$, we will call the
$f_c$ given above the \term{straigtened renormalization} of $f$. 

Suppose all periodic cycles of $f=f_c$ are repelling. Then if $f^n: U' \> U$ is a
(geometric) renormalization of $f$, then all of its periodic cycles are
repelling. It follows that all the periodic cycles of the map the
straightened renormalization $f_{c'}$ given by the preceeding theorem (so $h \circ
f^n = f_{c'} \circ h$ on $U'$) must also be repelling, because repelling periodic
cycles are preserved under quasiconformal conjugacy. Now, with the same supposition
on $f$, we will say that $f$ is $m$ $(m>1)$ times renormalizable if its straightened
renormalization is $m-1$ times renormalizable. (We say that $f$ is once
renormalizable if $f$ is renormalizable). We say that $f$ is \term{infinitely
renormalizable} if $f$ is $m$ times renormalizable for all $m$. If $f$ is not
infinitely renormalizable, then there exists a series of maps $f_0=f, f_1,\dots
f_m$ (all of the form $f(z) = z^2+c$) such that $f_{i+1}$ is the straightened
renormalization of $f_i$, and $f_m$ is not renormalizable. 

\newcrit{Need to say that $K_f=J_f$ in all the cases we're considering.}

\subsection{Renormalization and Holomorpic Removability}
Suppose that $f$ is renormalizable, with period $n$. Let $q$ be the number of
external rays that land at $\alpha$. We will consider two cases:
\begin{enumerate}
\item \label{prim} \term{Primitive} renormalization: $n > q$
\item \label{sat} \term{Satellite} renormalization: $n=q$
\end{enumerate}

In Case \ref{prim}, there exists a nondegenerate critical annulus $\critannulus
N$, just as in the non-re\-norm\-al\-izable case
\cite{Milnor:local:connectivity}, and 
furthermore we can find $M \ge N$ such that 
$f^n: \critpiece{M+n} \> \critpiece M$ is a non-trivial quadratic-like map
\cite{Milnor:local:connectivity}  (and therefore $J_{\R f}
\subset \critpiece k$ for all $k$).  In Case \ref{sat}, $\alpha \in \partial
\critpiece k$ for all $k$, and $\alpha \in J_{\R f}$. This will make Case \ref{sat}
a little harder to handle in what follows.

We will first prove:
\begin{prop} \label{pr1}
Suppose $f$ is primitively renormalizable, and $J_{\R f}$ is holomorphically
removable. Then $J_f$ is holomorphically removable. 
\end{prop}
\newcrit{Really, I'm talking about $K_f$. But the proof should actually work for
$J_f$.}
\begin{proof}
In this case, $\alpha \notin \closure{\critpiece {N+1}}$ (where $\critannulus N$
is non-degenerate, as mentioned above), so if $P$ is the piece of level $N$ such
that $\alpha \in \partial P$, and $P \subset \piece 0c$ (where
$c=f(0)$ is the critical value), then $c \notin P$, and then we can proceed
as in Subsection \ref{get slices}, and in fact the entire argument of Section
\ref{embedding rns} applies, so Lemma \ref{getS} applies, and therefore Lemma
\ref{nonuniform} applies, i.e, for all pieces $P$ of the Yoccoz puzzle for $f$,
$\dist {J\cap P}P < \infty$.

We can also prove the Tiling Lemma, \ref{tiling}, for $f$, as follows. 
 Let
$L=M+n$ (where $M$ is mentioned above). Then given $\critpiece p$, with $p>L$, we
let $R=J_{\R f}$, and let
$Q_i$ be the pieces of level $q_i>p$ such that $Q_i \subset \critannulus {q_i-1}$.
Now since $f^n:\critpiece{M+n} \> \critpiece M$ is a degree two branched cover, 
$f^n: \critannulus {r+n} \> \critannulus r$ is a degree two (unbranched) cover
for all $r \ge M$. Therefore, by induction, $f^{tn}: \critannulus {r+tn} \>
\critannulus r$ is an unbranched cover (of degree $2^t$) for all $r\ge M$. It
follows that $f^{L-q_i}$ is univalent on $q_i$, since we can find $tn$ such that
$M+n > q_i-tn \ge M$, and then $f^{tn} :\critannulus{q_i-1} \> \critannulus
{q_i-1-tn}$ is a covering, so $f^{tn}$ is univalent on any pieces $Q_i$ with $Q_i
\subset \critannulus{q_i-1}$. Letting $T=\critpiece P - {R \cup \bigcup
\closure{Q_i}}$ we find that $T \cap J = \nullset$ because $R =
\bigcap_{q>p}\critpiece q=\critpiece p-\bigcup_{q>p}\critannulus{q-1}$, and
$\critannulus {q-1}\cap J=(\bigcup_{q_i=q}Q_i) \cap J$. This completes the proof
of the tiling lemma. 

Now we can apply the proof of Lemma \ref{uniform} verbatim, and conclude that
there exists a $K$ such that $\dist {J\cap P}P \le K$ for all pieces $P$. We can
then apply the proof of Main Theorem,\ref{HR}, given in section \ref{udb}, but
there is one minor detail: the diameter of the pieces for $f$ do not go to zero.
However, given a homeomorphism $h: \cx \> \cx$ such that $h \on {\cx -J_f}$ is
conformal, we can still find a sequence of quasiconformal mappings $h_n: \cx \>
\cx$ such that $h_n=h$ on $\Gamma_n$ and on the unbouounded component of
$\cx-\Gamma_n$. Then, by the compactness of $K$-quasiconformal
mappings\cite{Ahlfors:book:qc, Lehto:Virtanen}, we can find a uniform limit $h_\infty$
of a subsequence of the $h_n$. Then $h_n=h$ on $\cx-J$ and also on 
$\{z \mid \exists n: f^n(z) =\alpha\}$, which is dense in $J$. So $h_\infty=h$ and
$h_\infty$ is $K$-quasiconformal, so we conclude that $h$ is always
$K$-quasiconformal, where $K$ depends only on $f$. Then we can conclude, as in
section \ref{udb}, that $J_f$ is holomorphically removable.
\end{proof}

\begin{cor}\label{pr2}
Suppose $f$ is finitely renormalizable (with all periodic cycles repelling), so
there exists a sequence $f_0=f, f_1, \dots, f_m$ where $f_{i+1}$ is the
straightened renormalization for $f_i$, and $f_m$ is non-renormalizable. Suppose
that each renormalization is primitive. Then $J_f$ is holomorphically removable. 
\end{cor}
\begin{proof}
We prove this by backwards induction on $i$. Certainly $J_{f_i}$ is holomorphically
removable if $i=m$. If $J_{f_{i+1}}$ is holomorphically removable, then $J_{\R
f_i}$ is too, because there is a quasiconformal map from one to the other. Then by
the proposition, $J_{f_i}$ is holomorphically removable. 
\end{proof}

\subsection{Satellite Renormaliztion}
\label{sat ren}
We must now consider the case of Case \ref{sat}, i.e. satellite
renormalization. In this case our goal is still to first prove piece-dependent
distortion bounds. We observe that a sufficient hypothesis for Lemma \ref{getS}
(cf. Subsection \ref{mapping rns}) is the following:
\begin{hypothesis} \label{Cantor set of ray-pairs}
There exists mappings $q_1, q_2: \mt \> S^1$ such that 
$\forall x \in \mt,q_1(x) \simeq q_2(x)$, and $q_1,q_2$ extend to quasisymmetric
mappings $q_1:[0,1]\>[A,B], q_2:[0,1]\>[C,D]$ (where $q_2$ is
orientation-reversing), where $CS(A,B,C,D)$  is a vertical slice, as in section
\ref{mapping rns}. 
\end{hypothesis}
We then prove the following lemma:
\bl
Suppose that $f$ is finitely renormalizable (with all periodic cycles
repelling). Then Hypothesis \ref{Cantor set of ray-pairs} is satisfied for $J_f$.
\el
\begin{proof}
If $f$ is primitively renormalizable, then, as discussed previously (in the
proof of \ref{pr1}), the hypothesis holds. If $f$ is satellite renormalizable, we
will need the following lemma:
\bl
Suppose $f$ has a satellite renormalization, and let $\hat{f}$ be the
straightened renormalization of $f$. Then $J_f$ satisfies Hypothesis \ref{Cantor set
of ray-pairs} if $J_{\hat f}$ does. 
\el
\begin{proof}
There are two things to prove:
\begin{enumerate}
\item If $J_g$ satifies Hypothesis \ref{Cantor set of ray-pairs}, then it also
does for a slice $CS(A',B',C',D')$ with $\beta \in \partial CS(A',B',C',D')$ (so
$A'=D'=0$). 
\item If $J_{\hat f}$ satisfies the above conclusion, then $J_f$ satisfies
Hypothesis \ref{Cantor set of ray-pairs}
\end{enumerate}

To prove the first, we first apply $g^{-1}$ to get the slice in $\piece
0\beta$($=\critpiece 0$). Then, we can keep applying the branch of $g^{-1}$ that fixes $\beta$ to
get a series of slices limiting on $\beta$. Then we can map the first half of
the Cantor set to the first slice, and the third quarter of $\mt$ to the second
slice, and the seventh eighth of $\mt$ to the third slice, and so forth, and
then map the last point of $\mt$ to $\beta$. 

To prove the second, we need a folk result, which relates the combinatorial
ray-pairs for $J_{\hat f}$ to those of $J_f$\cite{Milnor:hairiness}. It states
that there exists a pair $(a_0,a_1)$ of binary strings, such that if $t_1 \equiv
t_2$ on $J_{\hat f}$, then $E(t_1) \equiv E(t_2)$, where $E$ is defined by 
$E(.d_1 d_2 d_3 \dots)=.a_{d_1} a_{d_2} a_{d_3}$. From this, the second step
easily follows, after we note that $E(.000\dots)$ and $E(.111\dots)$ are both
rays landing at $\alpha$, in the case where $f$ is primitively renormalizable.
\end{proof}

The result then follows by induction on the number of times that $f$ is
renormalizable. 
\end{proof}

Now given this lemma, we observe that if $f$ is finitely renormalizable, then
Lemma \ref{getS} holds for the pieces of the Yoccoz puzzle for $f$, and hence so
does Lemma \ref{nonuniform}. Now, we wish to proceed as in the case of
primitive renormalization, assuming that $J_{\hat f}$ is holomorphically
removable, and proving that $J_f$ is. As before $J_{\R f}$ is qc equivalent to
$J_{\hat f}$ and is hence holomorphically removable, but we run into a minor
glitch in proving the Tiling Lemma, because now $J_{\R f} \not\subset \critpiece
k$ for any $k$, so $J_{\R f} \cap \critpiece k$ is not compact. We can get
around this by proving the following lemma:
\begin{lem}
Suppose that $A \subset \cx$ is compact and holomorphically removable, and $U
\subset \cx$ is open, and $\partial U$ is locally holomorphically removable, 
and $\dist {A\cap U} U < \infty$. Then if $h:\closure{U} \> \cx$
is an embedding, and $h \on {U-A}$ is conformal, then $h \on U$ is conformal.
\end{lem}
(A closed subset $B \subset \cx$ is \term{locally holomorphically removable} if,
for all open sets $V \subset \cx$, if $h: B \> \cx$ is an embedding, and $h \on
{V-B}$ is conformal, then $h \on B$ is conformal.)

\noindent 
\begin{proof}
Given $h$ as above, we can find $\tilde h:\closure U \> \cx$ such that $\tilde h
\on U$ is quasiconformal, and $\tilde h\on U = h \on U$. Then $\tilde h^{-1}
\circ h$ maps $\closure U$ homeomorphically to itself and is the identity on
$\partial U$;  extend it by the identity to a homeomorphism $g:\cx \> \cx$. Then
$g$ is quasiconformal on $\cx-{A \cup \partial U}$, and therefore is quasiconformal
on $\cx -A$, because $\partial U$ is locally holomorphically removable. Then $g$
is qc on all of $\cx$, so $h\on U=\tilde h \circ g$ is qc, and hence conformal
(since $A$ must have measure 0). 
\end{proof}

Note that the boundary of any Yoccoz puzzle piece is locally holomorphically
removable, since every point has a neighborhood that is a smooth curve, or the
union of a point and something locally holomorphically removable. \newcrit{Why
don't I just prove that $J_f$ is locally holomorphically removable. Then I won't
need to go through this ridiculous nonsense.} So then we can apply the above
lemma with $A=J_{\R f}$ and $U=\critpiece k$, and then the proof of the tiling lemma goes
through (with $R=J_{\R f} \cap \critpiece k$---it's okay that $R$ is not
compact, since it's still HR in $\critpiece k$, as in the above lemma), and then
we can proceed just as in the primitive case, to get the analog of Lemma
\ref{pr1} (and hence Lemma \ref{pr2}). This completes the proof of holomorphic
removability of Julia sets of finitely renormalizable quadratic polynomials.

\section{Conjectures on Holomorphic Removability}
\label{conj}
The techniques used here to show holomorphic removability could conceivably have
much wider application. Boundaries of John domains have been shown already to be
holomorphically removable\cite{Jones:hr}; it seems that these techniques could
provide a different, and in some ways more elementary, proof. A careful
examination of how distortion bounds are obtained for the canonical model
mentioned above suggest that such bounds could be shown for much more general
models, and then some very general criterion for holomorphic removability could
be described, perhaps in terms of the capacities of certain sets.  There are
also further possible dynamical applications. Certainly it should be possible to
apply these techniques to obtain holomorphic removability for all higher degree
polynomial Julia sets where the Yoccoz theory yields local connectivity. It
seems likely that holomorphic removability can also be shown for Julia sets of
certain infinitely renormalizable quadratic polynomials for which local
connectivity and area zero are known
\cite{Lyubich, Lyubich:student}.  Finally, I conjecture that the
boundary of $M$ is itself holomorphically removable, and that it can
be proved to be so by this technique of cutting neighborhoods of the set into
pieces and showing distortion bounds for those pieces \cite{Kahn:msri}.

\bibliographystyle{math} 
\bibliography{thesis}

\end{document}